\documentclass[11pt]{amsart}
\usepackage[utf8]{inputenc}
\usepackage[T1]{fontenc}
\usepackage{tikz-cd}
\usepackage{enumitem}
\usepackage{amsmath}
\usepackage{amssymb}
\usepackage{amsthm,amsfonts}
\usepackage{xspace}
\usepackage{graphicx}
\usepackage{hyperref}
\hypersetup{hidelinks}
\usepackage{color}
\usepackage[linesnumbered,algoruled,boxed,lined]{algorithm2e}
\usepackage[noend]{algpseudocode}
\usepackage{tikz}
\usetikzlibrary{positioning,decorations.markings}
\usepackage{todonotes}
\usepackage{comment}
\usepackage{mathtools}
\usepackage[capitalise]{cleveref}
\usepackage{cite}
\usepackage{bbm}
\usepackage[labelformat=simple]{subcaption}
\newtheorem*{rep@theorem}{\rep@title}
\newcommand{\newreptheorem}[2]{%
\newenvironment{rep#1}[1]{%
 \def\rep@title{#2 \ref{##1}}%
 \begin{rep@theorem}}%
 {\end{rep@theorem}}}

 \newtheorem*{repone@theorem}{\repone@title}
 \newcommand{\newreponetheorem}[2]{%
\newenvironment{repone#1}[1]{%
 \def\repone@title{#2 \ref{##1}\,(i)}%
 \begin{repone@theorem}}%
 {\end{repone@theorem}}}

\theoremstyle{plain}
\numberwithin{equation}{section}
\newtheorem{theorem}{Theorem}[section]
\usepackage{thmtools}
\newreptheorem{theorem}{Theorem}
\newreptheorem{proposition}{Proposition}
\newreptheorem{corollary}{Corollary}
\newreponetheorem{corollary}{Corollary}
\newtheorem{lemma}[theorem]{Lemma}
\newtheorem{proposition}[theorem]{Proposition}
\newtheorem{corollary}[theorem]{Corollary}

\theoremstyle{definition}
\newtheorem{definition}[theorem]{Definition}
\newtheorem{remark}[theorem]{Remark}
\newtheorem{remarks}[theorem]{Remarks}
\newtheorem{example}[theorem]{Example}

\newtheorem{notation}[theorem]{Notation}

\DeclareMathOperator{\im}{im}
\DeclareMathOperator{\colim}{colim}
\DeclareMathOperator{\Ob}{ob}
\DeclareMathOperator{\grank}{grank}

\newcommand{\Down}[1]{\operatorname{dn}(#1)}
\newcommand{\DownCL}[1]{\operatorname{dn}[#1]}
\newcommand{\Up}[1]{\operatorname{up}[#1]}
\newcommand{\kbb}{\mathbbm k}
\newcommand{\fun}{G}
\newcommand{\ifun}{F}

\newcommand{\N}{\mathbb{N}}

\newcommand*{\rom}[1]{\expandafter\@slowromancap\romannumeral #1@}

\newcommand{\B}{\mathcal{B}}

\newcommand{\I}{I}

\newcommand{\gr}{\mathrm{gr}}

\newcommand{\x}{\mathbf{x}}
\newcommand{\y}{\mathbf{y}}
\newcommand{\ZZ}{\mathsf{Z}}

\newcommand{\Z}{\mathbb{Z}}

\renewcommand{\I}{I}

\DeclareMathOperator{\rank}{rank}
\DeclareMathOperator{\thin}{thin}
\DeclareMathOperator{\supp}{Supp}

\newcommand{\e}{\mathbf{e}}

\renewcommand{\Vec}{\mathbf{Vec}}

\definecolor{darkred}{rgb}{1, 0.1, 0.3}
\definecolor{darkblue}{rgb}{0.1, 0.1, 1}

\title[Limit Computation Over Posets via Minimal Initial Functors]{Limit Computation Over Posets \\ via Minimal Initial Functors}

\author{Tamal K. Dey}
\address{Department of Computer Science, Purdue University}
\email{tamaldey@purdue.edu}
\thanks{}
\author{Michael Lesnick}
\address{Department of Mathematics and  Statistics, University at Albany -- SUNY}
\email{mlesnick@albany.edu}

\begin{document}

\begin{abstract}
It is well known that limits can be computed by restricting along an initial functor,
and that this often simplifies limit computation.
We systematically study the algorithmic implications of this idea for diagrams indexed by a poset.  We say an initial functor $\ifun\colon C\to D$ with $C$ small is \emph{minimal} if the sets of objects and morphisms of $C$ each have minimum cardinality, among the sources of all initial functors with target $D$.  For $Q$ a finite poset or $Q\subseteq \N^d$ an interval (i.e., a convex, connected subposet), we describe all minimal initial functors $\ifun\colon P\to Q$.  In particular, we show that such $\ifun$ is always a subposet inclusion, up to isomorphism.  We give efficient algorithms to compute a choice of minimal initial functor.  In the case that $Q\subseteq \N^d$ is an interval, we give asymptotically optimal bounds on $|\hom P|$, the number of relations in $P$ (including identities), in terms of the number $n$ of minima of $Q$: We show that $|\hom P|=\Theta(n)$ for $d\leq 3$, and $|\hom P|=\Theta(n^2)$ for $d>3$.
We apply these results to give new bounds on the cost of computing $\lim \fun$ for a functor $\fun \colon Q\to \Vec$ valued in vector spaces.  For $Q$ connected, we also give new bounds on  the cost of computing the \emph{generalized rank} of $\fun$  (i.e., the  rank of the induced map $\lim \fun\to \colim \fun$), which is of interest in topological data analysis.
\end{abstract}

\maketitle

\tableofcontents

\section{Introduction}\label{Sec:Intro}
Let $\Vec$ denote the category of vector spaces over a fixed field $\kbb$.  Given a category $D$ with finitely many morphisms and a pointwise finite-dimensional functor $\fun\colon D\to \Vec$, a standard equalizer formula \cite[Theorem V.2.2]{mac2013categories} expresses $\lim \fun$ as the kernel of a matrix, which can be computed via Gaussian elimination.  However, if $D$ has many objects and morphisms and the vector spaces of $\fun$ have large dimensions, then this matrix will be large and computing its kernel may be costly.

In general, one can compute the limit of a functor $\fun\colon D\to E$ by restricting along an \emph{initial functor} $\ifun\colon C\to D$; see \cref{Sec:Initial_Functors}.  If the category $C$ has far fewer objects and morphisms than $D$, then restricting along $\ifun$ can substantially simplify the problem of computing $\lim \fun$.  While this idea is well known, to our knowledge its algorithmic implications have not been systematically studied.

In this paper, we study these implications in the case that $D=Q$ is a poset and $E=\Vec$.  (More generally, the core ideas of our approach apply to any category $E$ where products and equalizers can be computed.) We develop theory and algorithms both for the case that $Q$ is a finite poset and the case that $Q$ is  an \emph{interval} in $\N^d$ (i.e., a connected, convex subposet; see \cref{Def:Interval}).  Here, $\N^d$ is given the product partial order, i.e., $(y_1,\ldots,y_d)\leq (z_1,\ldots,z_d)$ if and only if each $y_i\leq z_i$.  As we explain below, our motivation for considering intervals in $\N^d$ arises from computational questions in topological data analysis (TDA), and specifically in multiparameter persistent homology \cite{carlsson2009theory,botnan2023introduction}.

We develop a theory of \emph{minimal} initial functors $\ifun\colon P\to Q$, i.e., those for which the sets of objects and morphisms of $P$ both have minimum cardinality, among the  sources of all initial functors with target $Q$.  We give algorithms for computing such functors and apply them to the problem of limit computation.  We also apply our results to the closely related problem of \emph{generalized rank computation}, discussed below.

We now summarize our main results, deferring the precise statements to \cref{Sec:Main_Results}.
\begin{enumerate}
\item We present a structure theorem for minimal initial functors, \cref{theorem:MIS_main}, which establishes the existence of a minimal initial functor $\ifun\colon P\to Q$ and gives an explicit description of all such $\ifun$, assuming $Q$ is either finite or an interval in $\N^d$.  Specifically, the theorem says that, up to canonical isomorphism, any such $\ifun$ is the inclusion of a certain type of subposet $P\subseteq Q$, which we call an \emph{initial scaffold}.  While $Q$ can have multiple non-isomorphic initial scaffolds, all have the same underlying set of objects (but different sets of relations).
\item For $Q\subseteq \N^d$ an interval, \cref{Thm:Size_of_Initial} bounds $|\hom P|$, the number of relations of an initial scaffold $P\subseteq Q$ (including identities), in terms of the number $n$ of minima of $Q$, as follows:
\[ |\hom P|=
\begin{cases}
\Theta(n) &\textup{ for }d\leq 3,
\\ \Theta(n^2)&\textup { for }d>3.
\end{cases}
\]
The $d=3$ case is arguably the most interesting one.  To obtain these bounds, we establish a connection between initial scaffolds of intervals in $\N^d$ and the support of Betti numbers of monomial ideals.  We then apply a bound on the Betti numbers of monomial ideals, due to Bayer, Peeva, and Sturmfels \cite{bayer1998monomial}.
\item For $Q$ any finite poset, we give an algorithm to compute a choice of initial scaffold $P\subseteq Q$ from the Hasse diagram of $Q$.
We also give two specialized algorithms to compute $P$ when $Q\subseteq \N^d$ is an interval, one for the case $d\leq 3$ and one for arbitrary $d$.  \cref{Thm:Compute_Hull_Skeleton,Thm:Computing_Hull_Skeleton_Intervals} bound the complexity of each algorithm.  Our algorithms for the interval case do not require the Hasse diagram of $Q$ as input, but instead take a more parsimonious representation of $Q$, called the \emph{upset presentation}; see \cref{Def:Upset_Presentation}.  If the interval $Q\subseteq \N^d$ is finite, then we can instead specify $Q$ via its sets of minima and maxima.
\item We apply these results to derive new algorithms and bounds for computing the (co)limit of a functor $\fun\colon Q\to \Vec$; see \cref{Thm:Gen_Lim_Comp,Thm:Interval_Lim_Comp}.  The statements of these bounds depend on how $\fun$ is input to the algorithms.  For one concrete corollary of our bounds, consider a finite interval $Q\subseteq \N^d$ with a total of $n$ minima and maxima and a free presentation $F_1\to F_0$ of $\fun$ where $F_1$ and $F_0$ have maximum rank $r$.  Letting $\omega< 2.373$ denote the exponent of matrix multiplication, \cref{Thm:Chain_complex_limit,Rem:Drop_nr2_intro,rem:maxima_instead_of_upset_pres} imply that the cost of computing $\lim \fun$ is
\begin{alignat*}{2}
&O(n\log n+r^3) &\quad\textup{ for }d=2,\\
&O((nr)^{\omega}) &\quad\textup{ for }d=3,\\
&O(n^4+n^{\omega+1}r^{\omega}) &\quad\textup{ for } d>3.
\end{alignat*}
Our results on limit computation also encompass the case that $\fun$ is given implicitly as the homology of a chain complex of free functors, as well as the case that  matrix representations of the structure maps of $\fun$ are explicitly given.
    \item Given a connected, finite poset $Q$ and functor $\fun\colon Q\to \Vec$, we apply our results to bound the cost of computing the \emph{generalized rank} of $\fun$ \cite{kim2018generalized}, i.e., the rank of the induced map $\lim \fun\to \colim \fun$; see \cref{Thm:Grank_Gen_Hull_Skeleta,Thm:Grank_Nd_Hull_Skeleta}.   This invariant has a simple representation-theoretic interpretation: It is the number of \emph{interval summands} of $\fun$ with support $Q$ \cite[Corollary 3.2]{chambers2018persistent}.  As we discuss below, the problem of computing generalized ranks is of interest in TDA and has been studied in prior work \cite{DKM22,dey2024computing,asashiba2024interval}. Our bounds extend and improve on prior bounds by Dey, Kim, and Mémoli \cite{DKM22} and by Dey and Xin \cite{dey2024computing}.
\end{enumerate}

\subsection*{Motivation from Topological Data Analysis}
 Our motivation for considering (co)limits and generalized ranks over intervals in $\N^d$ arises from \emph{multiparameter persistence} \cite{carlsson2009theory,botnan2023introduction}, a subarea of TDA where functors $\fun\colon \N^d\to \Vec$ are the central algebraic objects of study.  For $d>1$, such functors have wild representation theory  \cite[Section 8]{botnan2023introduction}, making them relatively challenging to work with in practical applications.  This is in contrast to the case $d=1$, where each $\fun\colon \N\to \Vec$ decomposes in an essentially unique way into very simple indecomposables called \emph{interval modules} \cite{webb1985decomposition}, yielding a complete  invariant called  the \emph{persistence diagram} or \emph{barcode}.

To cope with the wildness of multiparameter persistence, the TDA community has sought to identify computable invariants of a functor $\fun\colon \N^d\to \Vec$ for use in data analysis.  Many invariants have been proposed \cite[Section 4.2]{botnan2023introduction}.  Such invariants can be used for data visualization \cite{lesnick2015interactive,botnan2024signed,vipond2018multiparameter} or as features for supervised machine learning \cite{vipond2018multiparameter,carriere2020multiparameter,corbet2019kernel,vipond2021multiparameter,beltramo2022euler,loiseaux2023framework,loiseaux2023stable,DS25,xin2023gril}.

One of the most widely studied invariants is the \emph{rank invariant} \cite{carlsson2009theory}, which is the map sending each pair $(a\leq b)\in \N^d\times \N^d$ to $\rank(G_{ab})$.  In the case $d=1$, this is equivalent to the barcode of $G$.  The rank invariant fails to distinguish pairs of $\N^2$-indexed functors with very simple structure \cite[Example 2.2]{lesnick2015interactive}; this motivates the search for a more discriminating invariant.  One natural choice is the \emph{generalized rank invariant}  \cite{kim2018generalized}, which maps each interval $i\colon Q\hookrightarrow \N^d$ to the generalized rank of $\fun \circ i$.

Recent work \cite{DS25,xin2023gril} uses the generalized rank invariant to vectorize 2-parameter persistent homology for supervised learning.  To control computational cost, the generalized rank is computed only along intervals of $\N^2$ with a bounded number of minima and maxima.  These computations are enabled by an efficient algorithm of Dey, Kim, and Mémoli~\cite{DKM22} for computing generalized ranks of intervals in $\N^2$.  Experimental results reported in \cite{DS25,xin2023gril} demonstrate that this approach can perform competitively in supervised learning applications.  To extend the approach to the setting $d=3$, one needs an efficient algorithm for computing generalized ranks of intervals in $\N^3$.  This was the initial motivation for our work.

\subsection*{Related Work}
Several recent articles have studied generalized rank computation.  As mentioned above, Dey, Kim, and Mémoli \cite{DKM22} give an efficient algorithm to compute generalized ranks over intervals in $\N^2$.   The approach reduces the problem to a zigzag persistence computation.  Several algorithms for such a computation have been given in the TDA literature  \cite{maria2014zigzag,milosavljevic2011zigzag,carlsson2010zigzag}, including a recent one of Dey and Hou \cite{DH22} that is notable for its superior practical performance.   While neither initial functors nor their duals, final functors, are explicitly considered in \cite{DKM22}, that paper's approach can be understood in our general framework of minimal initial and final functors; we develop this viewpoint in \cref{Sec:Main_Compute_Gen_Rank}.

 Inspired by the results of \cite{DKM22}, Dey and Xin \cite{dey2024computing} give an algorithm to compute the generalized rank over a general poset  by reducing the problem to a zigzag persistence computation. However, it turns out to be asymptotically more efficient to compute the generalized rank by directly computing the limits and colimits; compare \cref{Thm:Grank_Gen_Hull_Skeleta} and \cite[Theorem 6.2]{dey2024computing}.  In fact, this is true even if one does not exploit the speedups enabled by initial and final functors.

For $Z$ a finite poset, $i\colon Q \hookrightarrow Z$ an interval, and $\fun\colon Z\to \Vec$, Asashiba and Liu  \cite{asashiba2024interval} give an explicit formula for the generalized rank of $\fun\circ i$ as the difference of the ranks of two matrices.  In that work, the formalism of (co)limits is not explicitly used, and no explicit bound is given on the cost of computing the generalized rank.

A recent paper of Brüstle, Desrochers, and Leblanc \cite{brustle2025generalized} also considers the problem of generalized rank computation via initial and final functors, and contains results related to ours.  In particular,  \cite[Theorem C]{brustle2025generalized} is related to our structure theorem, \cref{theorem:MIS_main}, but concerns full subposets; see \cref{Sec:Brustle}.  In addition, \cite{brustle2025generalized} independently obtains a variant of one of our technical results, \cref{Lem:grank_initial_final}, which we use to bound the cost of computing the generalized rank over an interval in $\N^2$.  In contrast to our work, \cite{brustle2025generalized} does not explicitly deal with algorithmic questions, nor does it give specialized results for the case of  intervals in $\N^d$.

Applying a Möbius inversion to the generalized rank invariant of a finitely generated functor $\fun\colon \N^d\to \Vec$ yields an equivalent invariant called the \emph{generalized persistence diagram (GPD)}, introduced by Kim and Mémoli \cite{kim2018generalized}.  The GPD has been actively studied in recent work \cite{botnan2024signed,kim2024bigraded,kim2025super,carriere2025sparsification,asashiba2023approximation} as part of a larger effort to develop a theory of \emph{signed barcodes} for generalized persistent homology; see \cite[Section 9]{botnan2023introduction} for an overview.  While the GPD is generally too large to be computed in practice  \cite{kim2025super}, the construction clarifies how the generalized rank invariant relates to the usual 1-parameter persistence diagram.

\subsection*{Organization}
 \cref{Sec:Preliminaries} covers preliminaries.  \cref{Sec:Main_Results} is the backbone of the paper: Here, we precisely state the main definitions and results of our work, together with examples and context.   The remaining sections of the paper give the algorithms and proofs underlying our main results:  \cref{Sec:Structure_of_Minimal_Functors} gives the proof of \cref{theorem:MIS_main}, our structure theorem for minimal initial functors; \cref{Sec:Size_Bounds} gives the proof of our size bound \cref{Thm:Size_of_Initial}; \cref{Sec:Compute_Initial_Scaffold} presents our algorithms for computing initial scaffolds, and proves the bounds on their complexity; \cref{Sec:Limit_Comp_Proofs} completes the proofs of our bounds on limit computation; and \cref{Sec:Main_Compute_Gen_Rank} proves our bounds on generalized rank computation.
\cref{Sec:Future_Work} briefly discusses directions for future work.

\subsection*{Acknowledgments} We thank Ezra Miller for pointing out that \cref{Essential_and_Betti_1} follows from a version of Hochster's formula \cite[Theorem 1.34]{miller2005combinatorial}.  We also thank the anonymous reviewers for their insightful and very helpful feedback.
TD acknowledges the support of NSF grants DMS 2301360 and CCF 2437030. ML acknowledges the support of a grant from the Simons Foundation (Award ID 963845).

We, the authors, developed all ideas of this paper and wrote all of the text ourselves.  Claude and ChatGPT provided invaluable proofreading input, catching many typos and several minor errors.  Claude also made several helpful suggestions about the exposition.

\section{Preliminaries}\label{Sec:Preliminaries}
In this section, we review background on categories, posets, limits, initial functors, and $\Vec$-valued functors that we will need throughout the paper.  While we have aimed for an accessible treatment of the necessary category theory, we do assume familiarity with a few very basic notions, such as functors and natural transformations.   For background on category theory, see, e.g., the textbooks \cite{riehl2017category,mac2013categories}.

\subsection{Categories and Posets}

We denote the collections of objects and morphisms of a category $C$ as $\Ob C $ and $\hom C$, respectively.  For $c,c'\in \Ob C$, the set of morphisms from $c$ to $c'$ is denoted $\hom(c,c')$. We use subscripts to denote the action of a functor on objects and morphisms.  That is, for a functor $F\colon C\to D$ and $c\in \Ob C $ we let $F_c=F(c)$, and for a morphism $\gamma\colon  c\to c'$, we let $F_\gamma=F(\gamma)$.  We call $C$ \emph{thin} if $|\hom(c,c')|\leq 1$ for all $c,c'\in \Ob C $.  If $C$ is thin, then $\gamma\colon c\to c'$ is determined by its source and target, so we write $F_{cc'}=F(\gamma)$.

Recall that a functor $F\colon C\to D$ is called \emph{faithful (respectively, full)} if for each $c,c'\in \Ob C$, the induced map $\hom(c,c')\to \hom(F_c,F_{c'})$ is an injection (respectively, surjection).

\begin{definition}\label{Def:Embedding}
An \emph{embedding} is a faithful functor that is injective on objects.
\end{definition}

\begin{definition}
A poset is a pair $(Q,\leq)$, where $Q$ is a set and $\leq$ is a binary relation on $Q$, called the partial order, satisfying the following properties:
\begin{itemize}
\item $q\leq q$ for all $q\in Q$,
\item if $p\leq q$ and $q\leq p$, then $p=q$,
\item if $p\leq q$ and $q\leq r$, then $p\leq r$.
\end{itemize}
When $p\leq q$ and $p\ne q$, we write $p<q$.  We often abuse notation slightly and let $Q$ denote the poset $(Q,\leq)$.  By a further slight abuse, we refer to the elements of $\leq$ as \emph{relations}.
\end{definition}
We regard a poset $Q$ as a thin category, where $\Ob Q$ and $\hom Q $ are the sets of elements and relations of $Q$, respectively.  A poset $P$ is called a \emph{subposet} of $Q$ if $\Ob P\subseteq \Ob Q$ and $\hom P\subseteq \hom Q$.  In this case, we write $P\subseteq Q$.  By a slight abuse of terminology, we sometimes also refer to the inclusion map $P\hookrightarrow Q$ as a subposet.  We call $P$ a \emph{full subposet} if the inclusion functor $P\hookrightarrow Q$ is full, i.e., $p\leq p'$ in $P$ whenever $p,p'\in \Ob P$ and $p\leq p'$ in $Q$.

Intersections of posets are defined in the obvious way, i.e., by taking intersections of both the underlying sets and the sets of relations.  We define the union of posets analogously.  In general, a union of posets need not be a poset, as it may fail the transitivity property, but all unions considered in this paper satisfy transitivity and thus are posets.

We say a relation $p<q$ in $Q$ is a \emph{cover relation} if it does not factor as $p<r<q$ for any $r\in Q$.
The \emph{Hasse diagram} of a finite poset $Q$ is the directed graph $(\mathcal V,\mathcal E)$, where $\mathcal V=\Ob Q$ and $\mathcal E=\{(p,q)\mid p<q \textup{ is a cover relation}\}.$

\begin{definition}\label{Def:Connected_Category}
Let $C$ be a category.
\begin{itemize}
\item[(i)] For $c,c'\in \Ob C$, a \emph{path} from $c$ to $c'$ in $C$ is a sequence \[c=c_1,\ldots,c_k=c'\] in $\Ob C$ such that either $\hom(c_i,c_{i+1})$ or $\hom(c_{i+1},c_i)$ is non-empty for each $i\in \{1,\ldots,k-1\}$.
\item[(ii)]
$C$ is \emph{connected} if $C$ is non-empty (i.e., has at least one object) and there exists a path between each pair of objects of $C$.  Otherwise, $C$ is \emph{disconnected}.
\item[(iii)] A \emph{component} of $C$ is a maximal connected subcategory of $C$.
\end{itemize}
\end{definition}

Note that a finite poset $Q$ is connected if and only if the undirected graph underlying its Hasse diagram is a connected graph.  Note also that each object or morphism of $C$ belongs to exactly one component of $C$.
\begin{definition}\label{Def:Interval}
For a poset $Z$, a non-empty subset $Q\subseteq Z$ is an \emph{interval} if
\begin{itemize}
\item $Q$ is connected,
\item whenever $p\leq q\leq r$ with $q\in Z$ and $p,r\in Q$, we have $q\in Q$.
\end{itemize}
We regard an interval $Q\subseteq Z$ as a full subposet of $Z$.
\end{definition}

\begin{definition}\label{Def:Zigzag_Poset}
Let $\ZZ$ be the
full subposet of $\Z^2$ with
\[\Ob(\ZZ)=\bigcup_{z\in \Z} \{(z,-z),(z,1-z)\}.\]
A \emph{zigzag poset} is a poset isomorphic to a finite interval in $\ZZ$.
\end{definition}

\begin{example}\label{ex:Zigzag}
The poset $P$ with the following Hasse diagram is a zigzag poset.
\begin{center}
\begin{tikzcd}[ampersand replacement=\&,row sep=1.7ex,column sep=2.5ex,scale=.3,severy label/.append style={font=\small}]
\& v \arrow{r}  \& w                      \&                   \&                                                                                               \\
\&   \& x \arrow[]{r}\arrow[]{u}   \& y 	 					            \\
\&    \&                                                             \& z,\arrow[]{u}
\end{tikzcd}
\end{center}
\end{example}

\begin{definition}
Given a poset $Q$,
\begin{itemize}
\item[(i)] the \emph{(closed) downset} of $q\in Q$ is the set \[\DownCL{q,Q}\coloneqq\{p\in Q\mid p\leq q\},\]
\item[(ii)] the \emph{open downset} of $q\in Q$ is \[\Down{q,Q}\coloneqq \DownCL{q,Q}\setminus \{q\}= \{p\in Q\mid p<q\},\]
\item[(iii)] generalizing (i), the \emph{downset of $Q'\subseteq Q$} is the set \[\DownCL{Q',Q}\coloneqq\{p\in Q\mid p\leq q \textup{ for some } q\in Q'\}.\]
\item[(iv)] We say $P\subseteq Q$ is a \emph{downset}   if $P=\DownCL{Q',Q}$ for some $Q'\subseteq Q$.
\end{itemize}
Henceforth, we regard each of these sets as a full subposet of $Q$.  Dualizing, we also obtain corresponding definitions of \emph{upsets}, which we denote using ``$\mathrm{up}$'' rather than ``$\mathrm{dn}$''.  For example, the upset of $Q'\subseteq Q$ is the set \[\Up{Q',Q}\coloneqq \{p\in Q\mid p\geq q \textup{ for some } q\in Q'\}.\]
\end{definition}

We let $M_Q$ denote the set of minima (i.e., minimal elements) of a poset $Q$.  When $Q$ is clear from context, we sometimes omit it from the above notation, e.g.,
\begin{align*}
\DownCL{q}&\coloneqq\DownCL{q,Q},\\
M&\coloneqq M_Q.
\end{align*}

\begin{definition}
A poset $Q$ is \emph{well-founded} if each non-empty subset of $Q$ has a minimal element.
\end{definition}

If $Q$ is well-founded, then $Q$ contains no infinite chain $q_1>q_2>q_3 > \cdots$.  Assuming the axiom of choice, the converse is also true.

\subsection{Limits}\label{Sec:Limits}
Given a small category $C$, a category $D$, and a functor $\fun\colon C\to D$, a \emph{cone on $\fun$} consists of an object $d\in \Ob D$, together with a morphism $\delta_c\colon d\to \fun_c$ for each $c\in \Ob C$, called a \emph{cone map}, such that for all morphisms $\gamma\colon c\to c'$ in $\hom C$, we have $\delta_{c'}=\fun_{\gamma}\circ \delta_c$.  We write the cone as $(d,\delta)$.

\begin{example}\label{ex:cone}
If $P$ is the poset of \cref{ex:Zigzag}, then a cone on a functor $\fun\colon P\to D$ is a commutative diagram in $D$ of the following shape, extending $\fun$:
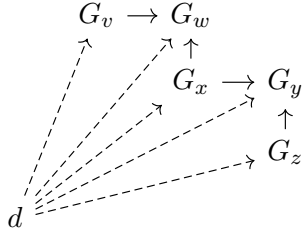
\begin{center}
\begin{tikzcd}[ampersand replacement=\&,row sep=1.7ex,column sep=2.5ex,scale=.3,severy label/.append style={font=\small}]
\&\fun_{v}\arrow{r}  \&\fun_{w}                      \&                   \&                                                                                               \\
\&   \&\fun_{x}\arrow[]{r}\arrow[]{u}   \&\fun_{y}	 					            \\
\&    \&                                                             \&\fun_{z}\arrow[]{u} \\
\textcolor{black}{d}\arrow[color=black,ruuu,dashed]\arrow[color=black,dashed,rruuu]\arrow[color=black,dashed,rruu]\arrow[color=black,dashed,rrruu]\arrow[color=black,dashed,rrru]
\end{tikzcd}
\end{center}
\end{example}

\begin{definition}
The \emph{limit} of $\fun$ is a cone $(d,\delta)$ such that for any other cone $(d',\delta')$, we have a unique morphism $\kappa\colon d'\to d$ satisfying $\delta'_c=\delta_c\circ\kappa$ for all $c\in \Ob C$.
\end{definition}
When they exist, limits are unique up to unique isomorphism, which justifies speaking of \emph{the} limit of $\fun$.  If $(d,\delta)=\lim \fun$, we sometimes abuse terminology slightly and refer to $d$ as $\lim \fun$, particularly when $\delta$ is clear from context.
Cocones and colimits are defined dually.

The category $\Vec$ of vector spaces over a fixed field is \emph{complete} and \emph{cocomplete}, meaning that $\lim \fun$ and $\colim \fun$ exist for any $\fun\colon C\to \Vec$.

Given functors $F\colon C\to D$ and $\fun\colon D\to E$ such that $\lim \fun$ and $\lim (\fun\circ F)$ both exist, the limit cone $(\lim \fun,\delta)$ on $\fun$ restricts to a cone $(\lim \fun,\delta')$ on $\fun\circ F$.  The universality of $\lim(\fun\circ F)$ then yields a morphism $\lim \fun\to \lim (\fun\circ F)$ such that the following triangle commutes for each $c\in \Ob C$, where the diagonal arrows are the limit cone maps:
\begin{center}
\begin{tikzcd} [scale=.3]
\lim \fun\arrow{r}\arrow{dr}              &\lim (\fun\circ F)     \arrow{d}       \\
        & \fun_{F_c}=(\fun\circ F)_{c}
\end{tikzcd}
\end{center}

This is functorial in the sense that, given another functor $F'\colon C'\to C$, the maps on limits form the following commutative triangle:
\begin{center}
\begin{tikzcd} [scale=.3]
\lim \fun\arrow{r}\arrow{dr}               &\lim (\fun\circ F)     \arrow{d}       \\
        & \lim(\fun \circ F\circ F')
\end{tikzcd}
\end{center}

\subsubsection{Limits as the Solutions to Linear Systems}\label{Sec:Lim_Equalizer}
For any small category $C$ and functor $\fun \colon C\to \Vec$, $\lim \fun$ admits a simple, concrete description as the solution of a system of linear equations: For $v\in \prod_{c\in \Ob C}\fun_c$ and $b\in \Ob C$, let $v_b$ be the projection of $v$ onto $\fun_b$.  For $\gamma\in \hom C$, write the source and target of $\gamma$ as $s(\gamma)$ and $t(\gamma)$, respectively.  A standard result \cite[Theorem V.2.2]{mac2013categories} is that
\begin{equation}\label{Eq:Concrete_Limit}
\lim \fun = \left\{v\in \prod_{c\in \Ob C} \fun_c\ : \  \fun_{\gamma}(v_{s(\gamma)})=v_{t(\gamma)}\ \forall\,  \gamma\in \hom C \right\},
\end{equation}
with each cone map $\lim \fun\to \fun_{b}$ the restriction of the projection \[\prod_{c\in \Ob C} \fun_c\twoheadrightarrow \fun_{b}.\]
Elements of the vector space on the right side of \cref{Eq:Concrete_Limit} are called \emph{sections of $\fun$}.

In the special case of a functor $\fun \colon Q\to \Vec$ with $Q$ a well-founded poset, we can express $\lim \fun$ as the solution to a smaller system of equations:
\begin{proposition}\label{Eq:Concrete_Limit_For_Finite_Poset}
For a well-founded poset $Q$ and $\fun\colon Q\to \Vec$, we have
\begin{equation}\label{Eq:Presections}
\lim \fun=\left\{v\in\! \prod_{m\in M} \fun_m\, :\, \fun_{lq}(v_l)=\fun_{mq}(v_{m})\ \forall\,  q\in Q,\   l,m\in (M\cap \Down{q})\right\},
\end{equation}
with each cone map $\lim \fun\to \fun_q$ given as the composition  \[\lim \fun \hookrightarrow \prod_{m\in M} \fun_m\twoheadrightarrow \fun_{l} \xrightarrow{\fun_{lq}} \fun_q,\] for any choice of $l\in M$ with $l\leq q$.
\end{proposition}

\begin{definition}\label{Def:presections}\mbox{}
\begin{itemize}
\item[(i)] A \emph{presection} of $\fun$ is an element of the vector space on the right side of \cref{Eq:Presections}.
\item[(ii)] A \emph{presection basis} of $\fun$ is a basis for this vector space.
\end{itemize}
\end{definition}

\begin{proof}[Sketch of proof of \cref{Eq:Concrete_Limit_For_Finite_Poset}]
Each presection of $\fun$ extends uniquely to a section of $\fun$ via the structure maps of $\fun$.  This defines an isomorphism from the vector space of presections to the vector space of sections.  It is immediate that this isomorphism commutes with the cone maps.
\end{proof}
\subsection{Initial Functors}\label{Sec:Initial_Functors}

\begin{definition}
Given a functor $F\colon C\to D$ and $d\in \Ob D$, the \emph{comma category} $(F\downarrow d)$ is the category whose
\begin{itemize}
\item objects are the pairs $(c,\gamma)$, where $c\in \Ob C$ and $\gamma\in \hom(F_c,d)$,
\item morphisms $(c,\gamma)\to (c',\gamma')$ are morphisms $\kappa \colon c\to c'$ in $C$ such that $\gamma=\gamma'\circ F_{\kappa}$.
\end{itemize}
\end{definition}

\begin{definition}\label{def:initial_functor}
A functor $F\colon C\to D$ with $C,D$ small is \emph{initial} if $(F\downarrow d)$ is connected for each $d\in \Ob D$.
\end{definition}
The following standard result says that we can compute a limit by restricting a diagram along an initial functor.

\begin{proposition}[{\cite[Theorem IX.3.1]{mac2013categories}}]\label{Prop:Initial_Functors_Preserve_Limits}
 Consider functors $F\colon C\to D$ and $\fun\colon D\to E$ with $C$ and $D$ small and $F$ initial.  If $\lim (\fun\circ F)$ exists, then $\lim \fun$ exists and the induced map $\lim \fun\to \lim (\fun\circ F)$ is an isomorphism in $E$.
\end{proposition}

\begin{remark}
The proof of \cref{Prop:Initial_Functors_Preserve_Limits} in \cite{mac2013categories} gives a concrete construction of the cone $\lim \fun$ from the cone $\lim(\fun\circ F)$, where each cone map is taken to be the composition of a cone map of $\lim(\fun\circ F)$ with a structure map of $\fun$.  This construction makes precise the idea that to compute $\lim \fun$, it suffices to compute $\lim(\fun \circ F)$.
\end{remark}

We will be particularly interested in functors $F\colon C\to Q$ where $Q$ is a poset.  In this case, for $q\in Q$, $(F\downarrow q)$ can be identified with the full subcategory of $C$ with objects $\{c\in \Ob C\mid F_c\leq q\}$.  In particular, this yields the following observation:

\begin{proposition}\label{Prop:Initial_Subposets}
An inclusion of posets $P\hookrightarrow Q$ is initial if and only if for each
$q\in Q$, the intersection poset  $\DownCL{q,Q}\cap P$ is connected.
\end{proposition}

\begin{example}\label{Ex:Multiple_Initial_Functors}
Let $Q$ be the poset with underlying set \[Q=\{t,u,v,w,x,y,z\}\] whose Hasse diagram  has edges
\[\{(t,x), (u,x),(u,y),(v,y),(x,z),(y,z), (w,z)\};\]
see \cref{fig:Hasse_Diagram_Ex}.
For each $i\in \{1,2,3\}$, let $P_i\subseteq Q$ be the subposet with $\Ob P_i=\Ob Q$, whose Hasse diagram has edges
\[\{(t,x), (u,x),(u,y),(v,y),(w,z),e_i\},\]
where $e_1=(t,z)$, $e_2=(u,z)$, and $e_3=(v,z)$; see \cref{fig:3_Initial_Scaffolds}.  Then each inclusion $P_i\hookrightarrow Q$ is an initial functor.
\tikzset{
    mid arrow/.style={
        postaction={
            decorate,
            decoration={
                markings,
                mark=at position #1 with {\arrow{stealth[width=7pt]}}
            }
        }
    },
    mid arrow/.default=0.6
}

\begin{figure}
\centering

\begin{subfigure}[t]{0.48\textwidth}
\centering
\begin{tikzpicture}[scale =.7, node distance=1.5cm,
     every node/.style={draw, circle,inner sep=1pt,
    minimum size=4.5mm},
     every path/.style={thick}]

    \node (m1) at (-2,0) {$t$};
    \node (m2) at (0,0) {$u$};
    \node (m3) at (2,0) {$v$};
    \node (m4) at (4,0) {$w$};
    \node (e1) at (-1,1.5) {$x$};
    \node (e2) at (1,1.5) {$y$};
    \node (e3) at (0,3) {$z$};

    \draw[mid arrow,shorten >=2pt, shorten <=2pt] (m1) -- (e1);
    \draw[mid arrow,shorten >=2pt, shorten <=2pt] (m2) -- (e1);
    \draw[mid arrow,shorten >=2pt, shorten <=2pt] (m2) -- (e2);
    \draw[mid arrow,shorten >=2pt, shorten <=2pt] (m3) -- (e2);
    \draw[mid arrow,shorten >=2pt, shorten <=2pt] (e1) -- (e3);
    \draw[mid arrow,shorten >=2pt, shorten <=2pt] (e2) -- (e3);
    \draw[mid arrow=.45,shorten >=2pt, shorten <=2pt] (m4) -- (e3);
\end{tikzpicture}
\caption{The Hasse diagram of the poset $Q$ from \cref{Ex:Multiple_Initial_Functors}.}
\label{fig:Hasse_Diagram_Ex}
\end{subfigure}
\hfill
\begin{subfigure}[t]{0.48\textwidth}
\centering
\begin{tikzpicture}[scale =.7, node distance=1.5cm,
     every node/.style={draw, circle,inner sep=1pt,minimum size=4.5mm},
     every path/.style={thick}]

    \node (m1) at (-2,0) {$t$};
    \node (m2) at (0,0) {$u$};
    \node (m3) at (2,0) {$v$};
    \node (m4) at (4,0) {$w$};
    \node (e1) at (-1,1.5) {$x$};
    \node (e2) at (1,1.5) {$y$};
    \node (e3) at (0,3) {$z$};

    \draw[mid arrow,shorten >=2pt, shorten <=2pt,blue] (m1) -- (e1);
    \draw[mid arrow,shorten >=2pt, shorten <=2pt,blue] (m2) -- (e1);
    \draw[mid arrow,shorten >=2pt, shorten <=2pt,blue] (m2) -- (e2);
    \draw[mid arrow,shorten >=2pt, shorten <=2pt,blue] (m3) -- (e2);
    \draw[mid arrow,shorten >=2pt, shorten <=2pt,gray!40!white] (e1) -- (e3);
    \draw[mid arrow,shorten >=2pt, shorten <=2pt,gray!40!white] (e2) -- (e3);
    \draw[mid arrow=.45,shorten >=2pt, shorten <=2pt,blue] (m4) -- (e3);

    \draw[mid arrow,shorten >=2pt, shorten <=2pt,red] (m1) .. controls (-2,1.9) .. (e3);
    \draw[mid arrow,shorten >=2pt, shorten <=2pt,red] (m2) -- (e3);
    \draw[mid arrow=.3,shorten >=2pt, shorten <=2pt,red] (m3) .. controls (1.75,1.55) .. (e3);
\end{tikzpicture}
\caption{The Hasse diagram of each subposet $P_i\subseteq Q$ has all of the blue edges and exactly one of the red edges.}
\label{fig:3_Initial_Scaffolds}
\end{subfigure}

\caption{}
\end{figure}
\end{example}

\subsection{Functors from Posets to Vector Spaces}
Recall that $\Vec$ denotes the category of vector spaces over a fixed field $\kbb$.  Given a poset $Q$, we often refer to a functor $\fun\colon Q \to \Vec$ as a \emph{$Q$-module}.
If $v \in \fun_q$, then we call $q$ the \emph{grade of $v$}, and write $q=\gr(v)$.
The $Q$-modules form an abelian category whose morphisms are the natural transformations.

A \emph{generating set} of a $Q$-module $\fun$ is a set $S \subseteq \bigsqcup_{q\in  Q}\fun_q$ such that for any $v \in \bigsqcup_{q\in Q}\fun_q$, we have
\[ v = \sum_{i=1}^k c_i \fun_{\gr(v_i)\gr(v)}(v_i)\] for some vectors $v_1, v_2, \dots, v_k \in S$ and scalars $c_1, c_2, \dots, c_k \in \kbb$.
The set $S$ is \emph{minimal} if for all $v\in S$, the set $S\setminus \{v\}$ does not generate $\fun$.  We say  $\fun$ is \emph{finitely generated} if there exists a finite generating set of $\fun$.  Clearly, if $\fun$ is finitely generated, then a minimal generating set of $\fun$ exists.

\begin{definition}\label{Def:Int_Mod}
For $Z$ a poset and $Q\subseteq Z$ an interval, define the \emph{interval module} $\kbb^{Q}\colon Z\to \Vec$ by
\begin{align*}
\kbb^{Q}_q &=
\begin{cases}
\kbb &\text{if } q\in Q, \\
0 &\text{otherwise,}
\end{cases}
& \kbb^{Q}_{pq}=
\begin{cases}
\mathrm{Id}_{\kbb} &{\text{if } p,q\in Q},\\
0 &{\text{otherwise}}.
\end{cases}
\end{align*}
\end{definition}
The interval module $\kbb^Q$ is easily checked to be indecomposable.

\subsubsection{Free Modules}\label{Sec:Free_Modules}
A $Q$-module $\fun$ is \emph{free} if there exists a multiset $\mathcal B$ of elements in $Q$ such that \[\fun\cong \bigoplus_{q\in \mathcal B}\, \kbb^{\Up{q,Q}}.\]
An elementary linear algebra argument shows that if such a multiset $\mathcal B$ exists, then it is unique; we call  $|\B|$ the \emph{rank of }$\fun$.
Define the function $\beta^\fun\colon Q\to \mathbb N\cup\{\infty\}$ by taking $\beta^\fun(q)$ to be the number of copies of $q$ in the multiset $\mathcal B$.
A \emph{basis of $\fun$} is a minimal generating set.  One can check that for any basis $B$ of $\fun$, \[\beta^\fun(q)=|\{s\in B\mid \gr(s)=q\}|.\]

Given a morphism  $\gamma\colon \fun\to \fun'$ of finitely generated free $Q$-modules and a choice of ordered bases $B=\{b_1,\ldots,b_n\}$ and $B'=\{b'_1,\ldots,b'_m\}$ for $\fun$ and $\fun'$, we  represent $\gamma$ as a matrix $[\gamma]$ with coefficients in $\kbb$, with each row and each column labeled by an element of $Q$, as follows:
\begin{itemize}
\item The $j^{\mathrm{th}}$ column of $[\gamma]$ is the unique one such that \[\gamma(b_j)= \sum_{i : \gr(b_i')\leq\gr(b_j)}  [\gamma]_{ij} \fun_{\gr(b_i')\gr(b_j)}(b_i')\]
and $[\gamma]_{ij}=0$ if $\gr(b_i')\not \leq\gr(b_j)$.
\item The $i^{\mathrm{th}}$ row label is $\gr(b'_i)$.
\item  The $j^{\mathrm{th}}$ column label is $\gr(b_j)$.
\end{itemize}
Regarding $\gamma$ as a functor $Q\times \{0<1\}\to \Vec$, one easily checks that $\gamma$ can be recovered up to natural isomorphism from $[\gamma]$.

\subsubsection{Free Resolutions and Betti Numbers}\label{Sec:Free_Resolution_Betti}
Given a poset $Q$,
a \emph{free resolution} $X$ of a $Q$-module $\fun$ is an exact sequence of free $Q$-modules \[X= \quad \cdots \xrightarrow{\partial_3}  X_{2}\xrightarrow{\partial_2} X_1\xrightarrow{\partial_1} X_0\] such that $\fun\cong \mathrm{coker}(\partial_1)$.
$X$ is \emph{minimal} if any free resolution of $\fun$ has a direct summand isomorphic to $X$.

By a standard structure theorem  \cite[Sections 1.7 and 1.9]{peeva2010graded}, if $\fun\colon \N^d \to \Vec$ is finitely generated, then there exists a minimal resolution $X$ of $\fun$ that is unique up to isomorphism, with each $X_i$ finitely generated.
Moreover, \emph{Hilbert's syzygy theorem} says that $X_i=0$ for $i>d$.  We denote the function $\beta^{X_i}\colon \N^d\to \N$ as $\beta_i^\fun$ and call its values the
 \emph{$i^{\mathrm{th}}$ (multigraded) Betti numbers of $\fun$}.

\section{Main Definitions and Results}\label{Sec:Main_Results}
In this section, we introduce the central definitions of this paper and state our main results, deferring the main proofs to subsequent sections.
\subsection{Structure of Minimal Initial Functors}\label{Sec:Main_Results_Min_Init_Fun}

\begin{definition}
An initial functor $\ifun\colon C\to D$ is \emph{minimal} if for any initial functor $F'\colon C'\to D$, we have
\[|\Ob C|\leq |\Ob C'| \quad \textup{and} \quad |\hom C|\leq |\hom C'|.\]
\end{definition}

For $Q$ a poset, let
\[I_Q=\{q\in Q\mid \Down{q} \textup{ is disconnected}\}.\]
Note that since we require a connected poset to be non-empty (\cref{Def:Connected_Category}\,(ii)), we have $M_Q\subseteq I_Q$.

 Our first main result, \cref{theorem:MIS_main} below, gives a simple characterization of minimal initial functors  $\ifun\colon P\to Q$, for a large class of posets $Q$.  Specifically, the theorem shows that, up to canonical isomorphism, such functors are the inclusions of certain subposets of $Q$ called \emph{initial scaffolds}, which we now define.

\begin{definition}\label{Def:Initial_Scaffold}
An \emph{initial scaffold} of a well-founded poset $Q$ is a subposet $P\subseteq Q$ of the following form:
\begin{itemize}
\item The set of elements of $P$ is $I_Q$.
\item For each $p\in P$ and each component $A$ of $\Down{p,Q}$, $P$ contains exactly one relation $m<p$ with $m\in A$. (Note that $m$ is a minimal element of $Q$.)
\end{itemize}
\end{definition}
By definition, all choices of an initial scaffold $P\subseteq Q$ have the same elements and the same number of relations.  However, the next example shows that $Q$ can have multiple non-isomorphic initial scaffolds.

\begin{example}\label{Ex:Multiple_Initial_Scaffolds}
In \cref{Ex:Multiple_Initial_Functors}, the  initial scaffolds of $Q$ are the posets $P_1$, $P_2$, and $P_3$.  Note that $P_1\cong P_3\not\cong P_2$.
\end{example}

\begin{example}\label{Rem:Zigzag_N^2}
Given an interval $Q\subseteq \N^2$, order $M_Q$ by $x$-coordinate, writing
$M_Q=(m_1,\ldots,m_k)$.  Then \[I_Q=M_Q\cup \{m_i\vee m_{i+1}\mid 1\leq i\leq k-1\},\] where $\vee$ denotes the join operator.  The poset $Q$ has a unique initial scaffold $P$, which is the full subposet of $\N^2$ with objects $I_Q$.  Note that $P$ is a zigzag poset; see \cref{Def:Zigzag_Poset}.  For instance, if $k=3$, then $P$ is isomorphic to the  poset of \cref{ex:Zigzag}.
\end{example}

\begin{theorem}[Structure of Minimal Initial Functors]\label{theorem:MIS_main} Let  $Q$ be a well-founded poset.
\begin{itemize}
\item[(i)] An initial scaffold $ P\subseteq Q$ exists.
\item[(ii)] For any initial scaffold $ P\subseteq Q$, the inclusion $P\hookrightarrow Q$ is a minimal initial functor.
\item[(iii)] If $Q$ has a finite initial scaffold and $F\colon C\to Q$ is a minimal initial functor, then $F$ is an embedding and $\im F$ is an initial scaffold of $Q$.
\end{itemize}
\end{theorem}

We prove \cref{theorem:MIS_main} in \cref{Sec:Structure_of_Minimal_Functors}.  \cref{Cor:Finite_Initial_Scaffold} below establishes that any interval $Q\subseteq \N^d$ has a finite initial scaffold.  Thus, if $Q$ is finite or $Q\subseteq \N^d$ is an interval, then $Q$ satisfies the conditions of \cref{theorem:MIS_main}.

 \begin{remark}\label{Rem:Minimality_Full_Subposet}
\cref{theorem:MIS_main}\,(ii) immediately implies a variant of the result for full subposets:  If $Q$ is a well-founded poset and $P$ is the full subposet of $Q$ with $\Ob P=I_Q$, then the inclusion $P\hookrightarrow Q$ is initial.  Moreover, for any full subposet $P'\subseteq Q$ whose inclusion into $Q$ is initial, we have $P\subseteq P'$.  As illustrated by \cref{Ex:Multiple_Initial_Scaffolds}, initial scaffolds needn't be full, so this is a strictly weaker minimality statement than \cref{theorem:MIS_main}\,(ii).
 \end{remark}

\begin{remark}[Related work of Brüstle et al.]\label{Sec:Brustle}
In October 2025, Brüstle et al. posted a paper giving a variant of \cref{Rem:Minimality_Full_Subposet} under a different condition on the poset $Q$ \cite[Theorem C]{brustle2025generalized}.  Previously, in October 2024, one of us, Lesnick, gave a talk on early versions of our results \cite{Lesnick2024InitialHulls}, which centered on the result of  \cref{Rem:Minimality_Full_Subposet} for posets with finite downsets.  As one of the authors of \cite{brustle2025generalized} was present for this talk, \cite{brustle2025generalized} attributes a version of their Theorem C to us.
\end{remark}

\begin{remark}
It is straightforward to check that for any well-founded poset $Q$ and initial scaffold $P\subseteq Q$, $P$ is the unique initial scaffold of itself.  This provides an intrinsic characterization of posets arising as initial scaffolds.
\end{remark}

\subsection{Size of Initial Scaffolds}\label{Sec:Main_Results_Size}
A natural way to study $|\hom P|$, the number of relations of $P$, is to compare it to the number $n$ of minima of $Q$.  The next example shows that $|\hom P|$ can be arbitrarily large, even when $n=2$.

\begin{example}\label{Ex:Unbounded_Size}
Let $Q=\{m_1,m_2,q_1,\ldots,q_k\}$, where $m_1$ and $m_2$ are minima, $q_1,\ldots,q_k$ are maxima, and $m_i\leq q_j$ for all $i$ and $j$.  Then $Q$ is the unique initial scaffold of itself.
\end{example}

However, in the special case that $Q\subseteq \N^d$ is an interval, we show that  $|\hom P|$ is controlled by $n$, as follows:

\begin{theorem}\label{Thm:Size_of_Initial}
For $Q$ an interval in $\N^d$ with $n$ minima and $P\subseteq Q$ an initial scaffold,
  \[|\hom P|=
  \begin{cases}
  \Theta(n)&\textup{ for }d\leq 3,\\
  \Theta(n^2)&\textup{ for }d>3.
  \end{cases}
  \]
\end{theorem}

We prove \cref{Thm:Size_of_Initial} in \cref{Sec:Size_Bounds}.  As noted in \cref{Sec:Intro}, we prove the upper bounds by applying an upper bound from \cite{bayer1998monomial} on the Betti numbers of monomial ideals.  For $d\leq 3$, the lower bound of \cref{Thm:Size_of_Initial} is essentially trivial, while for $d>3$ the lower bound is obtained by an explicit construction.

\begin{corollary}\label{Cor:Finite_Initial_Scaffold}
Any interval $Q\subseteq \N^d$ has a finite initial scaffold.
\end{corollary}

\cref{Cor:Finite_Initial_Scaffold} is immediate from \cref{Thm:Size_of_Initial} and the following standard result:

\begin{lemma}[Dickson's Lemma \cite{dickson1913finiteness}]\label{Lem:Dickson's}
Any full subposet of $\N^d$ has finitely many minima.
\end{lemma}

\begin{proof}
This follows from the fact that polynomial rings are Noetherian.  To elaborate, if $Q\subseteq \N^d$ is a full subposet, then $Q$ and $\Up{Q,\N^d}$ have the same minima.  The minima of $\Up{Q,\N^d}$ are exactly the grades of a minimal set of generators of an ideal in the ring of polynomials $\kbb[x_1,\ldots,x_d]$; see \cref{Prop:Interval_as_MI}.  Since this ring is Noetherian, every ideal is finitely generated.
\end{proof}

\subsection{Computing Initial Scaffolds}\label{Sec:Main_Results_Compute_Min_Init_Fun}
We next state our bounds on the complexity of computing an initial scaffold, deferring the proofs and the algorithms underlying the bounds to \cref{Sec:Compute_Initial_Scaffold}.
Our first bound is for an arbitrary finite poset:

\begin{theorem}\label{Thm:Compute_Hull_Skeleton}
Given the Hasse diagram $(\mathcal V,\mathcal E)$ of a  finite poset $Q$, we can compute an initial scaffold of $Q$ in time  $O(|\mathcal V||\mathcal E|)$.
\end{theorem}

 The algorithm underlying \cref{Thm:Compute_Hull_Skeleton} is a straightforward application of depth-first search.

In the case that $Q\subseteq \N^d$ is a finite interval with $d$ constant, every vertex is incident to $O(1)$ edges, so the bound of \cref{Thm:Compute_Hull_Skeleton} simplifies to $O(|\mathcal V|^2)=O(|\Ob Q|^2)$.  Since $|\Ob Q|$ may be quite large compared to the amount of data required to specify $Q$, one might hope to improve this bound.  To this end, we give specialized algorithms for computing an initial scaffold of an arbitrary interval $Q\subseteq \N^d$; these take as input a different (and usually smaller) representation of $Q$, which we now introduce, following  \cite[Definition 6.4.1]{miller2025homological}:

\begin{definition}\label{Def:Upset_Presentation}
Given an interval $Q\subseteq \N^d$, let $Q'=\Up{Q,\N^d}\setminus Q$.  Note that $Q'$ is an upset in $\N^d$.  We call the pair
$(M_Q,M_{Q'})$ of sets of minima the  \emph{upset presentation} of $Q$ and let $\|Q\|=|M_Q|+|M_{Q'}|$.
\end{definition}

An interval $Q\subseteq \N^d$ is determined by its upset presentation.
Informally, we think of $M_Q$ and $M_{Q'}$ as the sets of birth and death points of $Q$, respectively.  We emphasize that $\|Q\|$ and $|\hom Q|$ can be very different.  In particular, the former is always finite by \cref{Lem:Dickson's}, whereas the latter can be infinite.

\begin{example}\label{Ex:Canonical_Presentation}
We give the upset presentations of three intervals in $\N^2$:
\begin{itemize}
\item[(i)] If $Q=\{(0,0)\}$, then $M_{Q}=\{(0,0)\}$ and $M_{Q'}=\{(1,0),(0,1)\}$.
 \item[(ii)] If $Q=\N^2$, then $M_{Q}=\{(0,0)\}$ and $M_{Q'}=\emptyset$.
 \item[(iii)] If $Q=\{(z,0)\mid z\in \N\}$, then $M_{Q}=\{(0,0)\}$ and $M_{Q'}=\{(0,1)\}$.
 \end{itemize}
\end{example}

If $Q$ is a finite interval in $\N^d$, then one can instead specify $Q$ by its extrema (i.e., sets of minimal and maximal elements).
However, \cref{Ex:Canonical_Presentation}\,(ii) and (iii) illustrate that an infinite interval $Q\subseteq \N^2$ is not necessarily determined by its extrema.  Throughout, we assume that when a finite interval $Q\subseteq \N^d$ is specified by its extrema, the maxima and minima are given separately.

The following result, which we prove in \cref{Sec:Compute_Initial_Scaffold}, bounds the cost of computing an initial scaffold of an interval in $\N^d$.

\begin{theorem}\label{Thm:Computing_Hull_Skeleton_Intervals}
Given the upset presentation of an interval $Q\subseteq \N^d$ with $\|Q\|=n$, we can compute an initial scaffold of $Q$ in time
\begin{itemize}
\item[(i)] $O(n\log n)$ for $d=2,3$,
\item[(ii)] $O(n^4)$ for $d>3$.
\end{itemize}
In the case that $Q$ is finite, these bounds also hold when $Q$ is instead specified by its extrema, with $n$ the total number of these.
\end{theorem}

Separate algorithms underlie the bounds (i) and (ii) of \cref{Thm:Computing_Hull_Skeleton_Intervals}.  In the case $d=3$, the algorithm for (i) constructs the initial scaffold by iterating through two-dimensional slices of $\N^3$.  The main step in the algorithm for (ii) is to compute an initial scaffold of the subposet of $\Up{Q,\N^d}$ consisting of all minima and joins of pairs of minima.
It seems likely that for small enough $d$, say $d=4$, the bound of \cref{Thm:Computing_Hull_Skeleton_Intervals}\,(ii) could be improved by extending our algorithm for the case $d=3$, but this appears to not be trivial.

\subsection{Computing Limits via Initial Scaffolds}
We next consider the application of initial scaffolds to limit computation.  First, we specify the precise definition of limit computation that we use in this paper:

\begin{definition}[Limit Computation]\label{Def:Lim_Computation}
Let $Q$ be either a finite poset or an interval in $\N^d$, and let $\fun$ be a finitely generated $Q$-module.  In this paper, \emph{computing $\lim \fun$} means computing a presection basis of $\fun$ (\cref{Def:presections}).
\end{definition}

\begin{remarks}\label{Rems:Def_of_Lim_Computation}\mbox{}
\begin{itemize}
\item[(i)] While this definition of limit computation does not entail explicit computation of the cone maps of $\lim \fun$, \cref{Eq:Concrete_Limit_For_Finite_Poset} makes clear that the cone maps are readily obtained from what we compute: Given a presection basis of $\lim \fun$, each cone map to a minimal element is given by coordinate projection, while each cone map to a non-minimal element is given by composing a coordinate projection with a structure map of $\fun$.
\item[(ii)] Using \cref{Prop:Initial_Functors_Preserve_Limits}, it is easily checked that if $j\colon P\hookrightarrow Q$ is an initial scaffold of $Q$, then the presection bases of $\lim \fun$ and  $\lim (\fun\circ j)$ are identical.  Therefore, according to \cref{Def:Lim_Computation}, the problems of computing $\lim \fun$ and $\lim( \fun \circ j)$ are identical.
\end{itemize}
\end{remarks}

To give an algorithm for computing the limit of a functor $\fun\colon Q\to \Vec$, we must specify how such functors are represented in our computations.  Here is one useful way:
\begin{definition}\label{Def:Matrix_Rep}
Given a finite poset $Q$ and a pointwise finite-dimensional $Q$-module $\fun$, a \emph{matrix representation} of $\fun$ consists of:
\begin{enumerate}
\item a choice of ordered basis for each of the vector spaces $(\fun_q)_{q\in Q}$,
\item the matrix representation $[\fun_{pq}]$ of each structure map $\fun_{pq}$ of $\fun$, with respect to these bases.
\end{enumerate}
\end{definition}

To bound the complexity of limit computation, we use the following well-known linear algebra result:
\begin{proposition}[\cite{jeannerod2013rank}]\label{Omega_Linear_Solve}
 A $y\times z$ matrix $A$ of rank $k$ can be transformed to row echelon form in time $O(yzk^{\omega-2})$.  Hence, using the standard $O(z^2)$-time backsolve algorithm, a basis of $\ker(A)$ can be computed in time \[O(yzk^{\omega-2}+z^2).\]
\end{proposition}
We refer to the algorithm underlying \cref{Omega_Linear_Solve} as \emph{Gaussian elimination}, though it is different than the usual cubic-time version of Gaussian elimination.

Together, \cref{Omega_Linear_Solve} and \cref{Eq:Concrete_Limit_For_Finite_Poset} yield the following naive bound on the cost of limit computation, where \[r(\fun)\coloneqq \max_{q\in Q}\  \dim \fun_q.\]
\begin{proposition}\label{Prop:Limit_Computation_for_Finite_Posets}
Let $Q$ be a finite poset with $n$ minima.  Given a matrix representation of a $Q$-module $\fun$ with $r(\fun)=r$, $\lim \fun$ can be computed in time
\[O(|\hom Q|n^{\omega-1}r^{\omega}).\]
\end{proposition}

\begin{proof}
The solution to the system of equations in \cref{Eq:Concrete_Limit_For_Finite_Poset} is the kernel of a matrix with $O(r |\hom Q|)$ rows and $O(rn)$ columns.  (To see why the row bound holds, note that the linear equations of \cref{Eq:Presections} are indexed by triples $(q,l,m)$, where $q\in Q$ and $l,m\in M\cap \Down{q}$.  Since equality is transitive, it suffices to fix a choice of $l$ for each $q\in Q$, while letting $m$ vary.) Since the matrix has rank at most $rn$, the result follows from \cref{Omega_Linear_Solve}.
\end{proof}

We aim to refine \cref{Prop:Limit_Computation_for_Finite_Posets} by computing initial scaffolds.  To this end, given a poset $Q$ which is either finite or an interval in $\N^d$, as well as a finitely generated $Q$-module $\fun$, we compute $\lim \fun$ in three steps:

\begin{enumerate}
    \item compute an initial scaffold $j\colon P\hookrightarrow Q$ via \cref{Thm:Compute_Hull_Skeleton} or \cref{Thm:Computing_Hull_Skeleton_Intervals},
    \item compute a matrix representation of $\fun\circ j$,
    \item compute $\lim(\fun\circ j)$ by applying \cref{Prop:Limit_Computation_for_Finite_Posets}.
\end{enumerate}

Letting $\Gamma$ denote the cost of the second step, this approach leads to the following general bounds on the cost of limit computation:

\begin{corollary}\label{Thm:Gen_Lim_Comp}
Given the Hasse diagram $(\mathcal V,\mathcal E)$ of a finite poset $Q$ with $n$ minima and a $Q$-module $\fun$ with $r(\fun)=r$, we can compute $\lim \fun$ in time
\[O\left(|\mathcal V||\mathcal E|+|\hom P|n^{\omega-1}r^{\omega}\right)+\Gamma,\]
where $P$ is an initial scaffold of $Q$.
\end{corollary}

\begin{proof}
By \cref{Thm:Compute_Hull_Skeleton}, we can compute an initial scaffold $P\subseteq Q$ in time $O(|\mathcal V||\mathcal E|)$.  The result now follows from  \cref{Prop:Limit_Computation_for_Finite_Posets}.
\end{proof}

\begin{corollary}\label{Thm:Interval_Lim_Comp}
Given the upset presentation of an interval $Q\subseteq \N^d$ with $\|Q\|=n$ and a $Q$-module $\fun$ with  $r(\fun)=r$, we can compute $\lim \fun$ in time
\begin{itemize}
\item[(i)] $O((nr)^\omega)+\Gamma$ for $d=2,3$,
\item[(ii)] $O(n^4+n^{\omega+1}r^\omega)+\Gamma$ for $d>3$.
\end{itemize}
\end{corollary}

\begin{proof}
By \cref{Thm:Computing_Hull_Skeleton_Intervals},  we can compute an initial scaffold $P\subseteq Q$ in time $O(n\log n)$ for $d=2,3$ and $O(n^4)$ for $d>3$.  In addition, \cref{Thm:Size_of_Initial} tells us that $|\hom P|=\Theta(n)$ if $d=2,3$ and $|\hom P|=\Theta(n^2)$ if $d>3$.  The result now follows from  \cref{Prop:Limit_Computation_for_Finite_Posets}.
\end{proof}

  The cost $\Gamma$ depends on the format in which the input diagram $\fun$ is given.  For example, if we are given a matrix representation of $\fun$, then $\Gamma$ is negligible.
 However, in the TDA applications that initially motivated this work, we typically are not given a matrix representation of $\fun$; instead, we are given a free presentation of $\fun$ or, more generally, a chain complex of free functors whose homology is isomorphic to $\fun$.  In this case, we may apply \cref{Prop:Chain_Complex_Compute_Matrix_Rep} below.
\begin{definition}\label{Def:QrComplex}
For $Q$ a poset, a \emph{$(Q,r)$-complex} is a pair of labeled matrices $([f],[g])$ representing a chain complex
\[X\xrightarrow{f} Y\xrightarrow{g} Z\]
of free $Q$-modules of maximum rank $r$, as in \cref{Sec:Free_Modules}.
We call $\ker g/\im f$ the \emph{homology} of the $(Q,r)$-complex.
\end{definition}

\begin{proposition}\label{Prop:Chain_Complex_Compute_Matrix_Rep}
Suppose we are given
\begin{enumerate}
\item a finite poset $P$, specified by its sets of elements and relations,
\item a poset $Q$, represented in such a way that pairs of elements can be compared in constant time,
\item a functor $j \colon P\to Q$, specified by its values on objects, and
\item a $(Q,r)$-complex with homology $H$.
\end{enumerate}
Then we can compute a matrix representation of $H \circ j$ in time $O(|\hom P|r^{\omega})$.
\end{proposition}

We prove \cref{Prop:Chain_Complex_Compute_Matrix_Rep} in \cref{Sec:From_Chain_Complex_to_Matrix_Rep} via straightforward linear algebra.

\begin{remark}
Note that the $(Q,r)$-complex $([f],[g])$ of \cref{Prop:Chain_Complex_Compute_Matrix_Rep} is in fact a free presentation of $H$ whenever $g=0$.  In the general case where $g$ may be non-zero, one can in principle first compute a (minimal) presentation of $H$, and then take this as the input to the algorithm underlying \cref{Prop:Chain_Complex_Compute_Matrix_Rep}.

In the case $Q=\N^2$ and in the case of $0^{\mathrm{th}}$ simplicial homology of filtrations indexed by arbitrary posets, efficient algorithms are known for computing a minimal presentation of $H$ \cite{morozov2025computing,lesnick2022minimal,fugacci2023compression,bauer2023efficient}.  In the case $Q=\N^2$, it has been observed that the size of a minimal presentation of $H$ is often far smaller than the size of the original chain complex \cite{fugacci2023compression}.
For $Q=\N^d$, a minimal presentation of $H$ can be computed via Gr\"obner basis algorithms, e.g., Schreyer's algorithm and its variants \cite{schreyer1980berechnung,cox1998using,erocal2016refined,la1998strategies}, but existing open-source implementations of such algorithms are expensive on the types of large, sparse input considered in TDA  \cite{lesnick2022minimal}.
\end{remark}

\begin{remark}
While some TDA constructions yield a chain complex whose chain modules are not free \cite[Section 5]{botnan2023introduction}, such a chain complex can often readily be converted to a free chain complex with isomorphic homology \cite{lesnick2015interactive,chacholski2017combinatorial,bauer2025fast}.
\end{remark}

 \begin{remark}\label{Rem:Compute_all_relations}
Given the Hasse diagram $(\mathcal V,\mathcal E)$ of a finite poset $Q$, in time $O(|\mathcal V||\mathcal E|)$ we can compute a $|\mathcal V|\times |\mathcal V|$ binary matrix that explicitly represents the partial order on $Q$, by performing a breadth-first search starting from each $v\in \mathcal V$.
Using this matrix, we can check whether any relation $p\leq q$ holds in constant time, as required in \cref{Prop:Chain_Complex_Compute_Matrix_Rep}.
 \end{remark}

 Together,  \cref{Thm:Gen_Lim_Comp,Thm:Compute_Hull_Skeleton,Rem:Compute_all_relations,Prop:Chain_Complex_Compute_Matrix_Rep} imply the following:

\begin{corollary}\label{Thm:Chain_complex_limit_gen}
Given the Hasse diagram $(\mathcal V,\mathcal E)$ of a finite poset $Q$ with $n$ minima and a
$(Q,r)$-complex with homology $H$, we can compute $\lim H$ in time \[O\left(|\mathcal V||\mathcal E|+|\hom P|n^{\omega-1}r^{\omega}\right),\]
where $P$ is an initial scaffold of $Q$.
\end{corollary}
Similarly, \cref{Thm:Interval_Lim_Comp,Thm:Computing_Hull_Skeleton_Intervals,Prop:Chain_Complex_Compute_Matrix_Rep} imply items (ii) and (iii) of the following corollary, which is an analogue of \cref{Thm:Chain_complex_limit_gen} for intervals in $\N^d$.

\begin{corollary}\label{Thm:Chain_complex_limit}
Given an $(\N^d,r)$-complex with homology $H$ and the upset presentation of an interval $i\colon Q\hookrightarrow \N^d$ with $\|Q\|=n$, we can compute ${\lim  (H\circ i)}$ in time
\begin{itemize}
\item[(i)] $O(n\log n+r^3+nr^2)$ for $d=2$,
\item[(ii)] $O((nr)^{\omega})$ for $d=3$,
\item[(iii)] $O(n^4+n^{\omega+1}r^{\omega})$ for $d>3$.
\end{itemize}
\end{corollary}

We prove  \cref{Thm:Chain_complex_limit}\,(i) in \cref{Sec:Computing_Limit_Z^2} using a somewhat different approach than for parts (ii) and (iii): By \cref{Rem:Zigzag_N^2}, an interval in $\N^2$ has  a unique initial scaffold, which is a zigzag poset. This allows us to compute the limit using a zigzag persistence computation rather than the general limit computation of \cref{Prop:Limit_Computation_for_Finite_Posets}, and this turns out to be more efficient.

\begin{remark}\label{Rem:Drop_nr2_intro}
The term $nr^2$ in the bound of \cref{Thm:Chain_complex_limit}\,(i) can be eliminated, assuming the output of the computation is given in a simple compressed format, where the bases at multiple isomorphic vector spaces of $H\circ i$ are represented by a single basis.  See \cref{Rem:Drop_nr2_term} for details.
\end{remark}

\begin{remark}\label{rem:maxima_instead_of_upset_pres}
For $Q$ finite, \cref{Thm:Chain_complex_limit,Thm:Interval_Lim_Comp} both also hold if $Q$ is instead specified by its minima and maxima and $n$ is the total number of these.
\end{remark}

\begin{remark}\label{Rem:Dual_Colims_Final}
All of the above definitions and results dualize immediately, in particular, yielding a \emph{final scaffold}, which is a \emph{minimal final functor} that can be used to compute colimits.
\end{remark}

\subsection{Generalized Rank Computation}\label{Sec:Application_to_Generalized_Rank}
As an application of our main results, we consider the problem of computing the \emph{generalized rank} of a diagram of vector spaces, a problem arising in TDA that has previously been studied in \cite{DKM22,dey2024computing,asashiba2024interval}.

\begin{definition}[\cite{kim2018generalized}]\label{Def:Gen_Rank}
For a connected poset $Q$ and a functor $\fun\colon Q\to \Vec$, the composition of cone and cocone maps $\lim \fun\to \fun_q\to \colim \fun$ is independent of the choice of $q\in Q$.  The rank of this map, which we denote $\grank(\fun)$,  is called the \emph{generalized rank} of $\fun$.
\end{definition}

Dey, Kim, and Mémoli \cite{DKM22} gave the following bound on the cost of computing a generalized rank over an interval in $\N^2$:

\begin{theorem}[\cite{DKM22}]\label{Thm:GRank_N2}
Given an $(\N^2,r)$-complex with homology $H$ and the $n$ extrema of a finite interval $i\colon Q\hookrightarrow \N^2$, we can compute $\grank(H\circ i)$ in time $O((n+r)^{\omega})$.
\end{theorem}

By computing initial and final scaffolds, we obtain the following bounds, whose proofs we give in \cref{Sec:Main_Compute_Gen_Rank}:

\begin{corollary}\label{Thm:Grank_Gen_Hull_Skeleta}
Let $Q$ be a finite, connected poset with $n$ extrema and let $P^{\mathcal I}$, $P^{\mathcal F}$ be initial and final scaffolds of $Q$.  Given the Hasse diagram $(\mathcal V,\mathcal E)$ of $Q$ and  a $(Q,r)$-complex with homology $H$, we can compute $\grank(H)$ in time
\[O(|\mathcal V||\mathcal E|+(|\hom P^{\mathcal I}|+|\hom P^{\mathcal F}|)n^{\omega-1}r^\omega).\]
\end{corollary}

\begin{corollary}\label{Thm:Grank_Nd_Hull_Skeleta}
Given an $(\N^d,r)$-complex with homology $H$ and the $n$ extrema of a finite interval $i\colon Q\hookrightarrow \N^d$, we can compute $\grank(H\circ i)$ in time
\begin{itemize}
\item[(i)] $O(n\log n+r^{\omega})$ for $d=2$,
\item[(ii)] $O((nr)^{\omega})$ for $d=3$,
\item[(iii)] $O(n^4+n^{\omega+1}r^{\omega})$ for $d>3$.
\end{itemize}
\end{corollary}

Note that \cref{Thm:Grank_Nd_Hull_Skeleta}\,(i) slightly strengthens \cref{Thm:GRank_N2}.  As we explain in \cref{Sec:Main_Compute_Gen_Rank}, \cref{Thm:Grank_Gen_Hull_Skeleta} and \cref{Thm:Grank_Nd_Hull_Skeleta}\,(ii) and (iii) follow readily from our main results, via straightforward linear algebra.  To prove \cref{Thm:Grank_Nd_Hull_Skeleta}\,(i), we use an additional argument to cast the generalized rank computation as a zigzag persistence computation.  This argument is similar to the one in \cite{DKM22}, though we frame it in the general formalism of initial and final scaffolds.

\section{Proof of Structure Theorem for Minimal Initial Functors}\label{Sec:Structure_of_Minimal_Functors}
In this section, we prove \cref{theorem:MIS_main}.  We prepare for the proof with several definitions and preliminary results.

Let $Q$ be a well-founded poset.  In this section, all open and closed downsets are taken with respect to $Q$, so we suppress $Q$ in our notation for downsets.  Recall from \cref{Def:Initial_Scaffold} that for any initial scaffold $P\subseteq Q$, we have $\Ob P=I_Q$.  Recalling that $M_Q\subseteq I_Q$, let $E_Q=I_Q\setminus M_Q$, i.e.,
\[E_Q=\{q\in Q\mid \Down{q} \textup{ is disconnected and non-empty}\}.\]
We call elements of $E_Q$ \emph{essential}.

\begin{definition}
Given a functor $F\colon C\to D$, the \emph{image} of $F$ is the ordered pair
\[\im(F)\coloneqq (F(\Ob C ),F(\hom C)),\] where $F(\Ob C ) \subseteq \Ob D$ and $F(\hom C)\subseteq \hom D $.  We write the first and second components of the pair $\im(F)$ as $\Ob(\im F)$ and $\hom(\im F)$, respectively.
\end{definition}
We sometimes abuse notation slightly and, e.g., write $d\in \im F$ to mean $d\in \Ob(\im F)$, or $\gamma\in \im F$ to mean $\gamma\in \hom(\im F)$.

Note that $\im F$ is not necessarily a subcategory of $D$, since $F(\hom C)$ may not be closed under composition; see \cref{Ex:Image_Not_a_Cat} below.  We let $\overline{\im F}$ denote the subcategory of $D$ generated by $\im F$, i.e., the category obtained from $\im F$ by including the compositions of all finite chains of composable morphisms in $\im F$.

We omit the easy proof of the following result:

\begin{lemma}\label{Lem:Embedding_Subcategory}
If $F$ is an embedding (see \cref{Def:Embedding}), then $\im F$ is closed under composition of morphisms, i.e., $\overline{\im F}=\im F$.
\end{lemma}

Since the definition of a connected category (\cref{Def:Connected_Category}) does not involve composition of morphisms, it extends to any collections of objects and morphisms in an ambient category.  In particular, the notion of connectivity makes sense for $\im F$.  We use this to state the next proposition, whose straightforward proof we omit.
\begin{proposition}\label{Prop:Connected_Images}
For any functor $F\colon C\to D$ with $C$ connected, $\im F$ is also connected.
\end{proposition}

\begin{proposition}\label{Prop:Image_of_Initial_Functor_Contains_Hull_Skeleton}For a well-founded poset $Q$ and initial functor $F\colon C\to Q$,
\begin{itemize}
    \item[(i)] $I_Q\subseteq \Ob(\im F)$,
    \item[(ii)] there exists an initial scaffold $P\subseteq Q$ with $P\subseteq \overline{\im F}$.
\end{itemize}
\end{proposition}

\begin{proof}
To prove (i), first recall from \cref{Sec:Initial_Functors} that for any functor $F\colon C\to Q$ and $q\in Q$, $(F\downarrow q)$ can be identified with the full subcategory of $C$ with objects $\{c\in \Ob C\mid F(c)\leq q\}$.  If $F$ is initial, then by definition $(F\downarrow q)$ is connected, hence non-empty, so $M_Q\subseteq \im F$.

Next we show that $E_Q\subseteq \im F$.  To do so, we assume to the contrary that there exists $q\in E_Q$ with $q\not \in\im F$, and show that then $(F\downarrow q)$ is disconnected, a contradiction.
Let  ${\hat F\colon (F\downarrow q)\to Q}$ be the restriction of $F$ to $(F\downarrow q)$.
Note that $\im \hat F=\im F\cap \DownCL{q}$, where \[\im F\cap \DownCL{q}\coloneqq \left(\Ob(\im F)\cap \Ob (\DownCL{q}),\hom(\im F)\cap \hom (\DownCL{q})\right).\]
Since we assume $q\not \in \im F$, we have $\im F\cap \DownCL{q} \subseteq \Down{q}$.  Moreover, since $M_Q\subseteq \im F$ and  $Q$ is well-founded, $\im F$ has non-empty intersection with each component of $\Down{q}$.  Since $q$ is essential, $\Down{q}$ is disconnected.  It follows that $\im \hat F=\im F \cap \Down{q}$ is disconnected.  Thus, $(F\downarrow q)$ is disconnected by the contrapositive of \cref{Prop:Connected_Images}, a contradiction.  We conclude that $E_Q\subseteq \im F$, and hence that $I_Q\subseteq \Ob(\im F)$, establishing (i).

In view of (i), to prove (ii) it suffices to show that for each $q\in E_Q$ and component $A$ of $\Down{q}$, there exists a relation $m<q$ in
$\overline{\im F}$ with $m$ a minimal element of $A$; indeed, we can then construct the desired initial scaffold $P$ by choosing one such relation for each such $q$ and $A$.
Since $Q$ is well-founded and $M_Q\subseteq I_Q\subseteq \Ob(\im F)$ by (i), we have $A\cap \im F\ne \emptyset$.  Moreover, since $(F\downarrow q)$ is connected,  \cref{Prop:Connected_Images} implies that $\im \hat F=\im F\cap \DownCL{q}$ is connected.  Thus, since $q\in \im F$ by (i), there exists at least one relation $a_0<q$ in $\im F$ where $a_0\in A$.  If $a_0\in M_Q$, we are done.  Otherwise, by applying the same argument to any component of $\Down{a_0}$, we obtain a relation $a_1<a_0$ in $\im F$ with $a_1\in A$.  Since $Q$ is well-founded, by continuing in this way, we eventually obtain a chain of relations
\[m<\cdots <a_1< a_0<q\] in $\im F$ with $m\in M_Q\cap A$.  We then have $m<q$ in
$\overline{\im F}$ as desired.
\end{proof}

The following example shows that \cref{Prop:Image_of_Initial_Functor_Contains_Hull_Skeleton}\,(ii) is no longer true if we replace $\overline{\im F}$ with $\im F$ in the statement.
\begin{example}\label{Ex:Image_Not_a_Cat}
Consider the posets $R=\{a,b,c,d,e\}$ and $Q=\{w,x,y,z\}$ whose Hasse diagrams are as shown in \cref{fig:Image_Not_a_Poset}.
\begin{figure}
\begin{center}
\begin{tikzpicture}
    [every node/.style={draw, circle,inner sep=1.5pt,minimum size=4.5mm},
     every path/.style={thick},scale=.75]

    \node[draw=none] (R) at (-4,3.75) {\Large $R:$};
    \node (a) at (-2,3) {$a$};
    \node (b) at (0,4.5) {$b$};
    \node (c) at (0,3) {$c$};
    \node (d) at (2,4.5) {$d$};
    \node (e) at (4,3) {$e$};

    \draw[->,thick] (a) -- (b);
    \draw[->,thick] (c) -- (b);
    \draw[->,thick] (c) -- (d);
    \draw[->,thick] (e) -- (d);

    \node[draw=none] (Q) at (-4,0) {\Large $Q:$};
    \node (w) at (-2,-1.5) {$w$};
    \node (x) at (0,0) {$x$};
    \node (y) at (2,1.5) {$y$};
    \node (z) at (4,0) {$z$};

    \draw[->,thick] (w) -- (x);
    \draw[->,thick] (x) -- (y);
    \draw[->,thick] (z) -- (y);
    \end{tikzpicture}
    \end{center}
\caption{The initial functor $F\colon R\to Q$ of \cref{Ex:Image_Not_a_Cat} maps each element of $R$ to the element of $Q$ directly below it.  For example, $F(b)=F(c)=x$.}
\label{fig:Image_Not_a_Poset}
\end{figure}
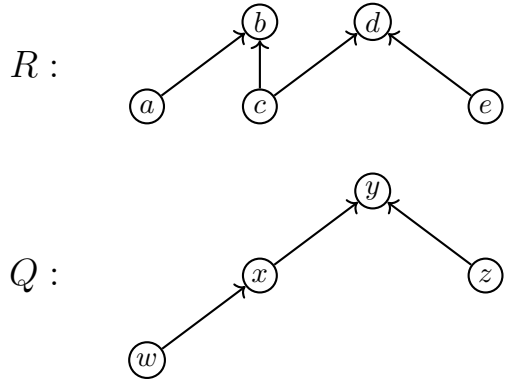
$Q$ has a unique initial scaffold $P\subseteq Q$, which is the full subposet with objects $\{w,y,z\}$.
Let $F\colon R\to Q$ be the functor given by \[F(a)=w,\  F(b)=F(c)=x,\ F(d)=y,\  F(e)=z.\]
 $F$ is easily checked to be initial, but $\im(F)$ does not contain the relation $w\leq y$, which is in $P$.
\end{example}

Note that for any category $C$, there is an associated thin category $\thin(C)$ with the same objects such that $\hom(c,c')$ is empty in $\thin(C)$ if and only if $\hom(c,c')$ is empty in $C$.  Moreover, we have a unique functor $\pi\colon C\to \thin(C)$ which is the identity on objects.  If $D$ is a thin category (e.g., a poset), then any functor $F\colon C\to D$ induces a functor $\thin(F)\colon \thin(C)\to D$ satisfying $F=\thin(F)\circ \pi$.

\begin{lemma}\label{Prop:Thinning_Initial}
If $F\colon C\to D$ is initial and $D$ is thin, then $\thin(F)$ is initial.
\end{lemma}

\begin{proof}
For any $d\in D$, the categories $(F\downarrow d)$ and $(\thin(F)\downarrow d)$ have the same objects because $C$ and $\thin(C)$ have the same objects, while $F$ and $\thin(F)$ act identically on objects.  In particular, since $(F\downarrow d)$ is non-empty, so is $(\thin(F)\downarrow d)$.  To show that $(\thin(F)\downarrow d)$ is connected, consider $(c,\gamma),(c',\gamma')$ in $(\thin(F)\downarrow d)$.  These are also objects in the connected category
$(F\downarrow d)$, so there is a path in $(F\downarrow d)$ connecting $(c,\gamma)$ and $(c',\gamma')$.  Applying $\pi$ yields a path in $(\thin(F)\downarrow d)$ connecting $(c,\gamma)$ and $(c',\gamma')$.  Hence $(\thin(F)\downarrow d)$ is connected.
\end{proof}

\begin{lemma}\label{Lem:Connectivity_And_Cover}
Given a poset $A$, let $\mathcal W$ be a set of connected subposets of $A$ with  $A=\bigcup_{P\in \mathcal W} P$.  Then $A$ is connected if and only if for each $x,y\in A$, there exists a sequence of posets $P_1,\ldots, P_l$ in $\mathcal W$ such that $x\in P_1$, $y\in P_l$, and  each $P_{i}\cap P_{i+1}$ is non-empty.
\end{lemma}

\begin{proof}
If $A$ is connected and $x,y\in A$, then there is a path \[x=a_1,\ldots,a_{l+1}=y\] in $A$.  By definition, $a_i$ and $a_{i+1}$ are comparable for each $i\in \{1,\ldots, l\}$.  Since $A=\bigcup_{P\in \mathcal W} P$, each pair $\{a_i,a_{i+1}\}$ belongs to some $P_i\in \mathcal W$.  The sequence $P_1,\ldots,P_{l}$ then has the desired properties.

Conversely, assume that for each  $x,y\in A$ there exists a sequence $P_1,\ldots,P_{l}$ of subposets in $\mathcal W$ as in the lemma statement.  Choose a sequence
\[x=a_1,a_2,\ldots,a_{l+1}=y\] such that  $a_{i+1}\in P_{i}\cap P_{i+1}$ for all $i\in \{1,\ldots, l-1\}$.  Then for each $i\in \{1,\ldots,l\}$, we have $a_{i},a_{i+1}\in P_i$, so since $P_i$ is connected, there exists a path in $P_i$ from $a_{i}$ to $a_{i+1}$.   These paths assemble into a path in $A$ from $x$ to $y$.  Hence $A$ is connected.
\end{proof}

We are now ready to prove \cref{theorem:MIS_main}.  We first recall the statement:

\begin{reptheorem}{theorem:MIS_main} Let  $Q$ be a well-founded poset.
\begin{itemize}
\item[(i)] An initial scaffold $ P\subseteq Q$ exists.
\item[(ii)] For any initial scaffold $ P\subseteq Q$, the inclusion $P\hookrightarrow Q$ is a minimal initial functor.
\item[(iii)] If $Q$ has a finite initial scaffold and $F\colon C\to Q$ is a minimal initial functor, then $F$ is an embedding and $\im F$ is an initial scaffold of $Q$.
\end{itemize}
\end{reptheorem}
\begin{proof}[Proof of {\cref{theorem:MIS_main}}]

For all $q\in Q$, each component $A$ of $\Down{q}$ contains a minimal element, since $Q$ is well-founded.   Hence an initial scaffold $P\subseteq Q$ always exists; this gives (i).  Alternatively, (i) follows from \cref{Prop:Image_of_Initial_Functor_Contains_Hull_Skeleton}\,(ii), since there always exists an initial functor with target $Q$, namely, the identity functor.

In view of
\cref{Prop:Initial_Subposets}, showing that the inclusion $P\hookrightarrow Q$ is initial amounts to checking that for each $q\in Q$, the intersection $\DownCL{q}\cap P$ is connected.
We proceed by induction on $q\in Q$ with respect to the partial order on $Q$; such an induction makes sense because $Q$ is well-founded.  The base case is that $q$ is minimal.  Then $q\in P$, so \[\DownCL{q}=\{q\}=\DownCL{q}\cap P\] is connected.
For the induction step, consider $q\in Q$ non-minimal and assume that $\DownCL{q'}\cap P$ is connected for all $q'<q$.
Let \[A_1,\ldots,A_k\] be the components of $\Down{q}$.

We claim that $A_i\cap P$ is connected for each $i$.   To show this, note that  $A_i$ is connected, as is each downset $\DownCL{a}$ for $a\in A_i$.
Thus, applying \cref{Lem:Connectivity_And_Cover} with $A=A_i$ and $\mathcal{W}=\{\DownCL{a}\}_{a\in A_i}$
gives that for each $x,y\in A_i\cap P$ there is a sequence $a_1,\ldots,a_l$ of elements in $A_i$ such that $x\in \DownCL{a_1}$, $y\in \DownCL{a_l}$, and
 \[\DownCL{a_i}\cap
\DownCL{a_{i+1}}\] is non-empty for each $i$.
Each such intersection contains a minimal element of $Q$ and hence an element of $P$, since $P$ contains all minima.
Thus, letting $P_i=\DownCL{a_i}\cap P$, each intersection $P_i\cap P_{i+1}$ is non-empty.
By the induction hypothesis, each $P_i$ is connected.   Noting that
\[A_i\cap P=\left(\bigcup_{a\in A_i} \DownCL{a}\right)\cap P=\bigcup_{a\in A_i} \DownCL{a}\cap P,\] applying \cref{Lem:Connectivity_And_Cover} with $A=A_i\cap P$ and $\mathcal W=\{\DownCL{a}\cap P\}_{a\in A_i}$  gives that $A_i\cap P$ is connected, as claimed.

If $q\not\in P$, then we have $k=1$ and so \[\DownCL{q}\cap P=\Down{q}\cap P=A_1 \cap P\] is connected.  If $q\in P$, then $q$ is essential.
Thus, by the definition of an initial scaffold,  for each $i$ we have a relation $m_i<q$ in $P$  with $m_i\in A_i$.  Moreover, $\DownCL{q}\cap P$ is the poset obtained from \[\Down{q}\cap P=\bigcup_{i=1}^k (A_i\cap P)\] by inserting $q$, together with the relations $m_i< q$ for each $i$.  Since each $A_i\cap P$ is connected, it follows that $\DownCL{q}\cap P$ is connected.  This finishes the proof that the inclusion $P\hookrightarrow Q$  is initial.

To see that $P\hookrightarrow Q$ is minimal, note that for any initial functor $F\colon C\to Q$, \cref{Prop:Image_of_Initial_Functor_Contains_Hull_Skeleton}\,(i) implies that $|\Ob P|=|I_Q|\leq |\Ob C|$.   Moreover, \cref{Prop:Image_of_Initial_Functor_Contains_Hull_Skeleton}\,(ii) implies that for every $q\in E_Q$ and component $A$ of $\Down{q}$, there is a relation $a<q$ in $\im F$ with $a\in A$.  It follows that $|\hom P |\leq |\hom C|$.  This proves (ii).

To prove (iii), assume that $Q$ has a finite initial scaffold (i.e., $I_Q$ is finite) and let $F\colon C\to Q$ be a minimal initial functor.   Then (ii) and \cref{Prop:Image_of_Initial_Functor_Contains_Hull_Skeleton}\,(i) together imply that $F$ is an injection on objects.  In addition, \cref{Prop:Thinning_Initial} implies that $C$ is thin, since otherwise we have $|\hom(\thin(C))|<|\hom C|$, contradicting that $F$ is minimal.  Therefore,  $F$ is faithful, hence an embedding.  Thus, we have $\overline{\im F}=\im F$ by \cref{Lem:Embedding_Subcategory}.  By \cref{Prop:Image_of_Initial_Functor_Contains_Hull_Skeleton}\,(ii), $\overline{\im F}=\im F$ contains an initial scaffold $P$.  The minimality of $F$ implies that $\im F=P$.
\end{proof}

In the induction underlying our proof of \cref{theorem:MIS_main}\,(ii), we have also essentially proven the following result, which will be useful in our algorithm to compute initial scaffolds of intervals in $\N^d$; see \cref{Sec:Compute_Hull_Skeleta_Nd}.

\begin{proposition}\label{Lem:Hull_Skeleton_Path_Components}
Given a well-founded poset $Q$, an initial scaffold $P\subseteq Q$, and $q\in Q$, the inclusion $P\cap \Down{q}\hookrightarrow \Down{q}$ induces a bijection on sets of components.
\end{proposition}

\begin{proof}
Since $Q$ is well-founded and $P$ contains all minima of $Q$, each component of $\Down{q}$ has non-empty intersection with $P$.  Thus, the map on components induced by the inclusion $P\cap \Down{q}\hookrightarrow \Down{q}$ is a surjection.  It is an injection if and only if $P\cap A_i$ is connected for each component $A_i$ of $\Down{q}$.  We have shown this in the inductive step of the proof of \cref{theorem:MIS_main}\,(ii).
\end{proof}

The following easy consequence of \cref{Lem:Hull_Skeleton_Path_Components} will be useful in our proof of \cref{Thm:Grank_Nd_Hull_Skeleta}\,(i) in \cref{Sec:Main_Compute_Gen_Rank}.

\begin{corollary}\label{Cor:Initial_Scaffold_Connected}
If $Q$ is a connected, well-founded poset, then any initial scaffold of $Q$ is connected.
\end{corollary}

\begin{proof}
Suppose $P$ is an initial scaffold of $Q$.  Let $Q'$ be the well-founded poset obtained from $Q$ by adding in a new maximum element $q$.  Then $q$ is not essential in $Q'$, since $\Down{q,Q'}=Q$ is connected.  Hence $P$ is also an initial scaffold of $Q'$.  \cref{Lem:Hull_Skeleton_Path_Components} implies that the inclusion $P=P\cap Q=P\cap \Down{q,Q'}\hookrightarrow \Down{q,Q'}=Q$ induces a bijection on sets of components.  Hence $P$ is connected.
\end{proof}

\section{Size of Initial Scaffolds of Intervals in {$\mathbb N^d$}}\label{Sec:Size_Bounds}
In this section, we prove \cref{Thm:Size_of_Initial}.  Recall that for an interval $Q\subseteq \mathbb \N^d$, this theorem gives tight asymptotic bounds on the number of relations of an initial scaffold of $Q$ in terms of the number of minima of $Q$.  The proof uses well-known results about the Betti numbers of monomial ideals.  We begin by reviewing these.

\subsection{Betti Numbers of Monomial Ideals}
\begin{definition}
Let $R^d=\kbb[x_1,\dots,x_d]$ denote the polynomial ring in $d\geq 1$ variables.  A \emph{monomial} is an element of $R^d$ of the form \[x^p\coloneqq x_1^{p_1}x_2^{p_2}\cdots x_d^{p_d}\] for some $p=(p_1,\ldots,p_d)\in \mathbb N^d$.  An ideal of $R^d$ generated by a set of monomials is called a \emph{monomial ideal}. \end{definition}

Given an upset $U\subseteq \N^d$, define the monomial ideal \[\hat U\coloneqq\langle x^p\mid p\in U\rangle\subseteq R^d.\]
We omit the straightforward proof of the following:
\begin{proposition}\label{Prop:Interval_as_MI}
The map $U\mapsto \hat U$ is a bijection between the set of upsets of $\mathbb N^d$ and the set of monomial ideals of $R^d$.
\end{proposition}

There is a standard isomorphism from the category of $\N^d$-modules to the category of $\N^d$-graded $R^d$-modules  \cite[Lemma 13]{corbet2018representation}.
By way of this isomorphism, we may regard a monomial ideal $J$ as an $\N^d$-module with
\[J_p=
\begin{cases} \kbb &\textup{ if } x^p\in J,\\
0 &\textup{ otherwise}
\end{cases}
\]
and $J_{pq}=\mathrm{Id}_{\kbb}$ whenever $x^p\in J$.  Thus, the definition of multigraded Betti numbers from \cref{Sec:Free_Resolution_Betti} applies to monomial ideals.

A version of \emph{Hochster's formula} \cite[Theorem 1.34]{miller2005combinatorial} gives the Betti numbers of a monomial ideal $J$ in terms of the reduced homology of a simplicial complex built from $J$.  To state the result, we need the following notation and definition: Given $s\subseteq \{1,\ldots,d\}$, let $\e_s=\sum_{j\in s} \e_j\in \{0,1\}^d$, where $\e_j\in \N^d$ denotes the $j^{\mathrm{th}}$ standard basis vector.

\begin{definition}
For $J\subseteq R^d$ a monomial ideal and $p\in \mathbb N^d$, define the \emph{upper Koszul complex of $J$ at $p$} to be the simplicial complex
\[K^p(J)=\{s\subseteq \{1,\ldots,d\} \mid s\ne \emptyset,\ x^{p-\e_s}\in J\}.\]
\end{definition}
Thus, vertices in $K^p(J)$ are indices $j\in \{1,\ldots,d\}$ such that $x^{p-\e_j} \in J$ and edges are pairs of indices $\{j,k\}$ such that $x^{p-\e_j-\e_k} \in J$. In the following, $\tilde H$ denotes reduced simplicial homology with coefficients in $\kbb$.
\begin{theorem}[{\cite[Theorem 1.34]{miller2005combinatorial}}]\label{Thm:Betti_Simplicial}
For any monomial ideal $J\subseteq R^d$, $p\in \N^d$, and $i\geq 1$,
\[\beta_i^J(p)=\dim \tilde H_{i-1}(K^p(J)).\]  Moreover, writing $J=\hat U$ for some upset $U\subseteq \N^d$,
\[\beta_0^J(p)=\begin{cases} 1 & \textup{ if $p$ is a minimal element of $U$},\\
0 & \textup{ otherwise}.
\end{cases}\]
\end{theorem}

\begin{corollary}
For any monomial ideal $J\subseteq R^d$ and $i\geq d$, we have $\beta_i^J=0$.
\end{corollary}

\begin{proof}
For all $i\geq d$ and $p\in \N^d$, the complex $K^p(J)$ has at most one $(i-1)$-simplex, which implies that $\tilde H_{i-1}(K^p(J))=0$.  The result then follows from \cref{Thm:Betti_Simplicial}.  Alternatively, one can obtain the result by applying Hilbert's syzygy theorem to the quotient module $R^d/J$.
\end{proof}

Let $\pi_0(W)$ denote the set of connected components of a simplicial complex $W$.  In the case $i=1$, standard properties of reduced homology yield the following corollary of \cref{Thm:Betti_Simplicial}:

\begin{corollary}\label{Cor:Betti_Components}
For any monomial ideal $J\subseteq R^d$ and $p\in \N^d$,

  \[\beta_1^J(p)=
 \begin{cases}
 |\pi_0(K^p(J))|-1 &\textup{ if }K^p(J)\ne \emptyset,\\
 0 &\textup{ if }K^p(J)=\emptyset.
 \end{cases}\]
\end{corollary}

\subsubsection{Upper Bounds on Betti Numbers}

A result of Bayer, Peeva, and Sturmfels \cite[Theorem 6.3]{bayer1998monomial} provides bounds on the Betti numbers of monomial ideals.  See also \cite{miller2005combinatorial} for an expository treatment.  We next state an asymptotic version of the result, leveraging a standard upper bound on the number of $i$-dimensional faces of a convex polytope \cite[Corollary 8.2]{edelsbrunner1987algorithms}.

In what follows, we assume $d$ is constant.

\begin{theorem}[\cite{bayer1998monomial,edelsbrunner1987algorithms}]\label{Thm:Betti_Bounds}
If $J\subseteq R^d$ is a monomial ideal with a minimal generating set of size $n$, then for each $i\in \{0,\ldots,d-1\}$, \[\sum_{p\in \mathbb N^d} \beta^{J}_i(p)= O(n^{\min(i+1,\lfloor \frac{d}{2}\rfloor)}).\]
\end{theorem}

\begin{proof}
 The argument of \cite[Theorem 6.3]{bayer1998monomial} gives that for $i\in \{0,\ldots,{d-1}\}$, the sum $\sum_{p\in \N^d} \beta^{J}_i(p)$ is at most the number of $i$-dimensional faces of some $d$-dimensional convex polytope  with $n$ vertices. According to \cite[Corollary 8.2]{edelsbrunner1987algorithms}, the number of such faces is $O(n^{\min(i+1,\lfloor \frac{d}{2}\rfloor)})$.
\end{proof}

\begin{remark}
 McMullen's famous \emph{upper bound theorem} \cite{mcmullen1970maximum} asserts that for all $n>d\geq 2$, the maximum number of $i$-dimensional faces of a $d$-dimensional convex polytope with $n$ vertices is attained for each $i$ simultaneously by a \emph{cyclic polytope}.  The bounds of \cite[Theorem 6.3]{bayer1998monomial} are formulated non-asymptotically, in terms of the number of $i$-dimensional faces of a cyclic polytope.  In fact, exact formulae for these numbers of faces are known \cite[Theorems 6.5 and 6.10]{edelsbrunner1987algorithms}.
\end{remark}

We do not know whether \cref{Thm:Betti_Bounds} is tight for $i\in \{2,\ldots,d-1\}$, but we observe that it is tight for $i=1$:

\begin{corollary}\label{Cor:Tight_Bounds_for_Betti_1}
For $J$ and $n$ as in \cref{Thm:Betti_Bounds}, we have
\[
\sum_{p\in \mathbb N^d} \beta^{J}_1(p)=
\begin{cases}
\Theta(n)& \textup{ for }d\leq 3,\\
\Theta(n^2)& \textup{ for }d> 3.
\end{cases}
\]
\end{corollary}

\begin{proof}
In the case $d=1$, we always have $n=1$ and $\beta_1^J=0$, so the bound $\Theta(n)$ holds trivially.  \cref{Thm:Betti_Bounds} implies that
\[
\sum_{p\in \mathbb N^d} \beta^{J}_1(p)=
\begin{cases}
O(n)& \textup{ for }d\in \{2,3\},\\
O(n^2)& \textup{ for }d> 3.
\end{cases}
\]
 It is easily checked that in the case $d=2$, we have \[\sum_{p\in \mathbb N^2} \beta^{J}_1(p)=n-1;\]
 see \cite[Proposition 3.1]{miller2005combinatorial} for a proof.
Thus, the bound $\Theta(n)$ holds for $d=2$.
 \cref{Ex:N4_tightness_of_betti_bound} below gives a family of monomial ideals in $R^4$ for which the sum of first Betti numbers grows quadratically with the size of a minimal generating set.
Thus, the bound $\Theta(n^2)$ holds for $d=4$.  Finally, for any monomial ideal $J\subseteq R^d$ and $d'>d$, note that $J'\coloneqq J\otimes_{\kbb} R^{d'-d}$ is an ideal of $R^d\otimes_{\kbb} R^{d'-d}=R^{d'}$ such that for all $i$,
\[\sum_{p\in \mathbb N^d} \beta^{J}_i(p)=\sum_{p\in \mathbb N^{d'}} \beta^{J'}_i(p).\]
Thus, the claimed bounds hold for all $d$.
\end{proof}

\subsection{Size Bounds}
For what follows, recall the definitions of the sets $I_Q$ and $E_Q$ given in \cref{Sec:Main_Results_Min_Init_Fun,Sec:Structure_of_Minimal_Functors}, respectively.
All upsets below are taken in the poset $\N^d$.

Our proof of \cref{Thm:Size_of_Initial} will use the following result, which we will also use in \cref{Sec:Compute_Hull_Skeleta_Nd,Sec:Computing_Hull_Skeleta_Of_Intervals_N3}, where we consider the computation of an initial scaffold of an interval in $\N^d$.

\begin{lemma}\label{Lem:Initial_Hulls_and_Upsets} For any interval $Q\subseteq \N^d$, we have
$\I_{Q}=\I_{\Up{Q}} \cap Q$.  In particular, $\I_{Q}\subseteq\I_{\Up{Q}}$.
\end{lemma}

\begin{proof}
It is clear that $M_Q=M_{\Up{Q}}=M_{\Up{Q}}\cap Q$.  Therefore, it suffices to check that $E_Q= E_{\Up{Q}}\cap Q$.

We claim that for all $q\in Q$, $\Down{q,Q}=\Down{q,\Up{Q}}$.  To check the claim, first note that $\Down{q,Q}\subseteq\Down{q,\Up{Q}}$ because $Q\subseteq \Up{Q}$.  Conversely, suppose $p\in \Down{q,\Up{Q}}$.  Since $\N^d$ is well-founded,  there exists $m\in M_{\Up{Q}}=M_{Q}$ with $m\leq p<q$.
Since $Q$ is an interval, we have $p\in Q$.  Thus, $\Down{q,\Up{Q}}\subseteq \Down{q,Q}$ which establishes the claim.

 The claim implies that for $q\in Q$, we have $q\in E_Q$ if and only if $q\in E_{\Up{Q}}$.
Thus, $E_Q=E_{\Up{Q}}\cap Q$, as desired.
\end{proof}

It can be shown that the containment $\I_{Q}\subseteq \I_{\Up{Q}}$ of \cref{Lem:Initial_Hulls_and_Upsets} is in fact an equality when $d\leq 2$, though we do not need this.  The next example shows that this containment  can be strict when $d=3$.
\begin{example}
Consider the following interval in $\mathbb N^3$,
\[Q=\{(1,0,0),(0,1,0),(0,0,1),(0,1,1),(1,0,1)\}.\]
We have \[\Up{Q}=\mathbb N^3\setminus\{(0,0,0)\}.\]  Note that $(1,1,0)$ is essential in $\Up{Q}$, so $(1,1,0)\in I_{\Up{Q}}$ but $(1,1,0)\not\in Q$. \end{example}

We denote the set of components of a small category $C$ as $\pi_0(C)$.

\begin{lemma}\label{Lem:Connectivity_of_G_and_Down_Z}
For $U\subseteq \mathbb N^d$ an upset and $p\in U$, let $T$ denote the full subposet of $\Down{p,U}$ with
\[\Ob T = \{y\in \Down{p,U} \mid y=p-\e_s \text{ for }s\subseteq \{1,\ldots,d\} \textup{ with }|s|\leq 2\}.\]
The inclusion $T\hookrightarrow \Down{p,U}$ induces a bijection $\pi_0(T)\to \pi_0( \Down{p,U})$.
\end{lemma}

\begin{proof}
In what follows, all downsets are taken in $U$.
Since $U$ is an upset, each component of $\Down{p}$    contains an element of the form $p-\e_i$ for some $i$.  Since $p-\e_i\in T$, the inclusion $T\hookrightarrow \Down{p}$ induces a surjection on components.  It remains to show that this inclusion is an injection.

Consider $w,y\in T$ lying in the same component of $\Down{p}$.  We need to show that $w,y$ lie in the same component of $T$.  There exists a path \[w=q_1,\ldots,q_k=y\] from $w$ to $y$ in $\Down{p}$.  By the transitivity property of a poset, we may assume that for each  triple $q_{i-1},q_i,q_{i+1}$ of consecutive points in the path, either \[q_{i-1}< q_i> q_{i+1}\quad \textup{or} \quad q_{i-1}> q_i< q_{i+1}.\]
We call $q_i$ an \emph{upper point} in the former case and a \emph{lower point} in the latter case.  By inserting an additional point into the path, if necessary, we may further assume that if $k>2$ and $q_1=p-\e_s$ for some $s$ with $|s|=2$, then $q_2$ is an upper point.  Symmetrically, we may also assume that if $k>2$ and $q_k=p-\e_s$
for some $s$ with $|s|=2$, then $q_{k-1}$ is an upper point.  Since $U$ is an upset, by replacing each upper point $q_i$ with a point above it, if necessary, we may take each upper point to be of the form $p-\e_j$.  Then for each lower point $q_i$, we have $q_{i-1}=p-\e_j$ and $q_{i+1}=p-\e_{j'}$ for some $j,j'$, so
\[\wedge(q_{i-1},q_{i+1})=
\begin{cases}
p-\e_j-\e_{j'} &\textup{ if }j\ne j',\\
p-\e_j &\textup{ if }j =j',
\end{cases}
\]
where $\wedge$ denotes the meet operator. Since $ p>\wedge(q_{i-1},q_{i+1})\geq q_i\in U$ and $U$ is an upset, we have $\wedge(q_{i-1},q_{i+1})\in \Down{p}$, hence $\wedge(q_{i-1},q_{i+1})\in T$.
Therefore, by replacing each lower point $q_i$ with  $\wedge(q_{i-1},q_{i+1})$, we obtain a path from $w$ to $y$ in $T$. This shows that $w,y$ lie in the same component of $T$,  completing the proof.
\end{proof}

\begin{lemma}\label{Prop:Downset_Components_And_Upper_Koszul_Components}
For $U\subseteq \mathbb N^d$ an upset and $p\in U$, there exists a bijection $\pi_0(\Down{p,U})\to \pi_0(K^p(\hat U))$.
\end{lemma}
\begin{proof}
\cref{Lem:Connectivity_of_G_and_Down_Z} provides a bijection $\pi_0(\Down{p,U})\to \pi_0(T)$, where $T$ is the subposet
in the statement of \cref{Lem:Connectivity_of_G_and_Down_Z}.  Note that $\pi_0(T)$ is in bijection with $\pi_0(G)$, where $G$ is the undirected graph underlying the Hasse diagram of $T$.  Let $K^+$ denote the barycentric subdivision of the 1-skeleton of $K^p(\hat U)$.  The simplicial map $f\colon G\to K^+$ given by $f(p-\e_s)=s$ is a graph isomorphism, and therefore induces a bijection $\pi_0(G)\to \pi_0(K^+)\cong \pi_0(K^p(\hat U))$.  Composing these bijections yields a bijection $\pi_0(\Down{p,U})\to \pi_0(K^p(\hat U))$.
\end{proof}

We denote the support of an $\N$-valued function $f$ as $\supp f$.

\begin{lemma}\label{Essential_and_Betti_1}\mbox{}
\begin{itemize}
\item[(i)] For $U\subseteq \mathbb N^d$ an upset and $p\in \N^d$,
\[
\beta^{\hat U}_1(p)=
\begin{cases}
|\pi_0(\Down{p,U})|-1
&\textup{ if } p\in U\textup{ and }\Down{p,U}\ne \emptyset,\\
0 &\textup{ otherwise}.
\end{cases}\]
\item[(ii)] We have $E_U= \supp \beta^{\hat U}_1$.
\end{itemize}
\end{lemma}

\begin{proof}
Item (i) follows immediately from \cref{Cor:Betti_Components} and \cref{Prop:Downset_Components_And_Upper_Koszul_Components}, as $K^p(\hat U)=\emptyset$ whenever $p\not \in U$.  Since $p\in E_U$ if and only if $\Down{p,U}$ has at least two components, (ii) follows from (i).
\end{proof}

\begin{proposition}\label{Initial_hull_of_Lattice_Interval}
For $U\subseteq \N^d$ an upset, we have $I_U=\supp \beta^{\hat U}_0\cup  \supp \beta^{\hat U}_1$.
\end{proposition}

\begin{proof}
\cref{Thm:Betti_Simplicial} implies that $M_U= \supp \beta^{\hat U}_0$, and by \cref{Essential_and_Betti_1}\,(ii), we have $E_U=\supp \beta^{\hat U}_1$.  Since $I_U=M_U\cup E_U$, the result follows.
\end{proof}

\begin{example}\label{Ex:N4_tightness_of_betti_bound}
We identify a family of upsets $(U^k)_{k\in \N}$ in $\N^4$ such that $U^k$ has $\Theta(k)$ minima and $\Theta(k^2)$ essential points.  By \cref{Essential_and_Betti_1}\,(ii), this corresponds to a family of ideals in $R^4$ for which the sum of first Betti numbers grows quadratically with the size of a minimal generating set, as needed in the proof of \cref{Cor:Tight_Bounds_for_Betti_1}.

Fixing $k\in \N$, let
\begin{align*}
A&=\{(i,k-i,0,0)\mid i\in \{0,\ldots,k\}\},\\
B&=\{(0,0,i,k-i)\mid i\in \{0,\ldots,k\}\},\\
U^k&=\Up{A\cup B}.
\end{align*}
To simplify notation, we let $U=U^k$.

Note that $M_U=A\sqcup B$, so $|M_U|=2(k+1)$.  We observe that for each $i,j\in\{0,\dots,k\}$, \[p_{ij}\coloneqq(i,k-i,j,k-j)\] is essential.  Indeed, writing
\begin{align*}
D_1=\{(i,k-i,a,b)\in \N^4\mid  (a,b)<(j,k-j)\},\\
D_2=\{(a,b,j,k-j)\in \N^4 \mid  (a,b)<(i,k-i)\},
\end{align*}
every point in $D_1$ is incomparable to every point in $D_2$. Thus, $\Down{p_{ij},U} = D_1\sqcup D_2$, which implies that $p_{ij}$ is essential.  Therefore, $|E_U|\geq (k+1)^2$, so by \cref{Thm:Betti_Bounds} and \cref{Essential_and_Betti_1}\,(ii), we have $|E_{U}|=\Theta(k^2)$, as claimed.
\end{example}

We are now ready to prove \cref{Thm:Size_of_Initial}.  Recall the statement:

\begin{reptheorem}{Thm:Size_of_Initial}
For $Q$ an interval in $\N^d$ with $n$ minima and $P\subseteq Q$ an initial scaffold,
  \[|\hom P|=
  \begin{cases}
  \Theta(n)&\textup{ for }d\leq 3,\\
  \Theta(n^2)&\textup{ for }d>3.
  \end{cases}
  \]
\end{reptheorem}

\begin{proof}
Suppose $P\subseteq Q$ is an initial scaffold of an interval $Q\subseteq \N^d$.  Note that for each $q\in Q$, $\Down{q,Q}$ has at most $d$ components.  Thus, for each $p\in P$, there are at most $d$ relations in $P$ of the form $m<p$.  Hence, $|\hom P|=O(|\Ob P|)$.  Therefore, to prove the result, it suffices to show that $|\Ob P|=|I_Q|$ satisfies the bounds of the theorem.  By \cref{Lem:Initial_Hulls_and_Upsets}, we have $|I_{Q}|\leq |I_{\Up{Q}}|$.  Since $\Up{Q}$ has the same minima as $Q$, \cref{Cor:Tight_Bounds_for_Betti_1} and \cref{Initial_hull_of_Lattice_Interval} together imply that
\[|I_{ Q}|=
\begin{cases}
\Theta(n)& \textup{ for }d\leq 3,\\
\Theta(n^2)& \textup{ for }d> 3.
\end{cases}
\]
The result follows.
\end{proof}

\section{Computing Initial Scaffolds}\label{Sec:Compute_Initial_Scaffold}

\subsection{Algorithm for Finite Posets}\label{Sec:Compute_Hull_Skeleton}
We next prove \cref{Thm:Compute_Hull_Skeleton} by giving an algorithm to compute the initial scaffold of a finite poset.  Recall the theorem statement:
\begin{reptheorem}{Thm:Compute_Hull_Skeleton}
Given the Hasse diagram $(\mathcal V,\mathcal E)$ of a  finite poset $Q$, we can compute an initial scaffold of $Q$ in time  $O(|\mathcal V||\mathcal E|)$.
\end{reptheorem}

\begin{proof}
The components of $Q$ can be computed in time $O(|\mathcal V|+|\mathcal E|)=O(|\mathcal V||\mathcal E|)$ by depth-first search on the Hasse diagram.  Given an initial scaffold for each component of $Q$, the union of these is an initial scaffold of $Q$.  It therefore suffices to treat the case of a connected poset $Q$.

Assuming $Q$ is connected, our algorithm to compute an initial scaffold $P\subseteq Q$ proceeds as follows:
We refine the partial order on $Q$ to a total order via the standard topological sorting algorithm in time $O(|\mathcal V|+|\mathcal E|)=O(|\mathcal E|)$.   We process the elements of $Q$ in increasing order.
When we visit $q\in Q$, we use depth-first search to compute the components of $\Down{q,Q}$, as well as a choice of minimal element $m_A\in A$ for each component $A$ of $\Down{q,Q}$; this takes time $O(|\mathcal E|)$.  If the number of components is not 1, then we add $q$ to $P$ and we also add the relation $m_A<q$ to $P$ for each $A$.  It is clear that the total time cost of this algorithm is $O(|\mathcal V||\mathcal E|)$.
\end{proof}

\subsection{Algorithm for Intervals in $\N^3$}\label{Sec:Computing_Hull_Skeleta_Of_Intervals_N3}
The rest of this section is devoted to the proof of \cref{Thm:Computing_Hull_Skeleton_Intervals}.  Let us recall the statement:

\begin{reptheorem}{Thm:Computing_Hull_Skeleton_Intervals}
Given the upset presentation of an interval $Q\subseteq \N^d$ with $\|Q\|=n$, we can compute an initial scaffold of $Q$ in time
\begin{itemize}
\item[(i)] $O(n\log n)$ for $d=2,3$,
\item[(ii)] $O(n^4)$ for $d>3$.
\end{itemize}
In the case that $Q$ is finite, these bounds also hold when $Q$ is instead specified by its extrema, with $n$ the total number of these.
\end{reptheorem}

Our approach is based on the following observation:

\begin{proposition}\label{Prop:Hull_Skeleton_as_Subset}
If $Q\subseteq \N^d$ is an interval and $P$ is an initial scaffold of $\Up{Q,\N^d}$, then  $P\cap Q$ is an initial scaffold of $Q$.
\end{proposition}

\begin{proof}
This is immediate from \cref{Lem:Initial_Hulls_and_Upsets} and the definition of initial scaffold.
\end{proof}

We prove \cref{Thm:Computing_Hull_Skeleton_Intervals}\,(i) by giving an algorithm that computes an initial scaffold of an interval $Q\subseteq \N^3$ in time $O(n\log n)$.  The algorithm and its analysis specialize to intervals $Q\subseteq \N^2$ by embedding $\N^2$ into $\N^3$ via the map $(x,y)\mapsto (x,y,0)$.

Letting $U=\Up{Q,\N^3}$, our algorithm first computes an initial scaffold $P$ of $U$ and then computes $P\cap Q$.  The computation of $P$ proceeds by iterating through planar slices of $U$.  To elaborate, for $z\in \N$, let
\[U_{z}=\{(x,y)\in \N^{2} \mid (x,y,z)\in U\},\]
and let $U_{-1}=\emptyset$.  Note that each $U_z$ is an upset, and that $U_{z}\subseteq U_{z+1}$ for all $z$.  To compute $P$, we  iterate through values of $z\in \N$ in increasing order, maintaining the sets of minima $M_{U_{z}}$.  When we pass from level $z-1$ to level $z$, we update $M_{U_{z-1}}$ to obtain $M_{U_{z}}$ and record all points $p\in E_U$ with $p_3=z$.  For each such $p$, we also compute the relations $m<p$ to be added to $P$.  To carry out these computations, we leverage a concrete levelwise description of $P$, given in \cref{Prop:Upset_to_Interval} below.

The computation of $P\cap Q$ from $P$ is done in a similar but much simpler way, also by iterating through planar slices.  The details of this depend on whether $Q$ is specified by its upset presentation or by its extrema, but are similar in the two cases.

\begin{notation}\label{Not:WzXz}\mbox{}
\begin{itemize}
\item[(i)] For $m=(m_1,m_2)\in \N^2$ and $t\in \N$, we identify $(m,t)$ with $(m_1,m_2,t)\in \N^3$.
\item[(ii)] For $m\in M_{U_z}$, let $m_{\downarrow}=(m,t)$, where $t\in \N$ is the unique value such that $(m,t)\in M_U$.
\item[(iii)] We order $M_{U_{z}}$ with respect to $x$-coordinate, writing
\[M_{U_z}=(m^1,\ldots,m^k).\]
\item[(iv)] Let \[X^z=\{(m,z) \mid m\in M_{U_{z-1}}\setminus M_{U_z}\}.\]
\item[(v)] Let  $V^z\subseteq \N^2$ be the set of joins of consecutive elements of $M_{U_{z}}$, i.e.,
\[V^z=\left\{m^i\vee m^{i+1}\mid i\in \{1,\ldots,k-1\}\right\}.\]
\item[(vi)] Let \[W^z= \{(m,z)\mid m\in V^z\setminus V^{z-1}\}.\]

\end{itemize}
\end{notation}
Note that since $M_{U_{-1}}=\emptyset$, we always have $V^{-1}=X^0=\emptyset$ and $W^0=V^0\times \{0\}$.

\begin{proposition}\label{Prop:Upset_to_Interval} \mbox{}
\begin{itemize}
\item[(i)] $E_U=\bigcup_{z\in \N} (X^z \cup W^z)$.
\item[(ii)] The relations of an initial scaffold $P\subseteq U$ are given as follows:
\begin{enumerate}
    \item If $p=(m,z)\in X^z$, then $P$ contains the relation $m_{\downarrow}<p$.
    \item If $p=(m,z)\in X^z\setminus W^z$, then $P$ also contains exactly one relation $(l,z)<p$, where $(l,z)$ is an arbitrary point of $M_{U}$ with $l<m$.
        \item If $p=(m^i\vee m^{i+1},z)\in W^z$,
    then $P$ contains the relations $m^i_{\downarrow} <p$ and $m^{i+1}_{\downarrow}<p$.
\end{enumerate}
\end{itemize}
\end{proposition}
We refer to the relations from (1)-(3) in \cref{Prop:Upset_to_Interval}\,(ii) as  \emph{vertical $X^z$-relations}, \emph{horizontal $X^z$-relations}, and \emph{$W^z$-relations}, respectively.

Before proving \cref{Prop:Upset_to_Interval}, we illustrate the result with an example.

\begin{example}\label{Ex:WzXz}
Let
\[
U=\Up{\{(0,6,0),(1,5,0),(3,4,0),(4,2,0),(5,0,0),(1,3,1),(2,2,1),(4,1,1)\},\N^3}.
\]
Then, as illustrated in \cref{fig:WzXz}, we have
\begin{align*}
M_{U_0}&=\left((0,6),(1,5),(3,4),(4,2),(5,0)\right),\\
M_{U_1}&=\left((0,6),(1,3),(2,2),(4,1),(5,0)\right),\\
X^1&=\{(1,5,1),(3,4,1),(4,2,1)\},\\
V^0&=\{(1,6),(3,5),(4,4),(5,2)\},\\
V^{1}&=\{(1,6),(2,3),(4,2),(5,1)\},\\
W^0&=\{(1,6,0),(3,5,0),(4,4,0),(5,2,0)\},\\
W^1&=\{(2,3,1),(4,2,1),(5,1,1)\}.
\end{align*}
Note that
\begin{align*}
X^1\cap W^1&=\{(4,2,1)\},\\
X^1\setminus W^1&=\{(1,5,1),(3,4,1)\}.
\end{align*}
Also note that $(1,6,1)\not\in W^1$ because $(1,6)\in V^0\cap V^1$.
\begin{figure}
\begin{center}
\begin{tikzpicture}[scale=0.7]
  \newcommand{\opendiamond}[3][0.30]{%
    \draw[black, thick]
      (#2,#3+#1) -- (#2+#1,#3) -- (#2,#3-#1) -- (#2-#1,#3) -- cycle;
  }

    \newcommand{\opensquare}[3][0.30]{%
    \draw[black, thick]
      (#2+#1,#3+#1) -- (#2+#1,#3-#1) -- (#2-#1,#3-#1) -- (#2-#1,#3+#1) -- cycle;
  }

  \foreach \x in {1,...,5} {
    \foreach \y in {1,...,6} {
      \fill[gray!40] (\x,\y) circle (1.2pt);
    }
  }

  \draw[->] (-0.2,0) -- (5.5,0) node[right] {$x$};
  \draw[->] (0,-0.2) -- (0,6.5) node[above] {$y$};

  \foreach \x in {1,...,5} {
    \draw (\x,.08) -- (\x,-0.08);
    \node[below, font=\scriptsize] at (\x,-0.08) {\x};
  }

  \foreach \y in {1,...,6} {
    \draw (.08,\y) -- (-0.08,\y);
    \node[left, font=\scriptsize] at (-0.08,\y) {\y};
  }

  \fill[blue] (5,0) circle (3pt);
  \fill[blue] (4,2) circle (3pt);
  \fill[blue] (3,4) circle (3pt);
  \fill[blue] (1,5) circle (3pt);
  \fill[blue] (0,6) circle (3pt);

  \fill[red] (1,3) circle (3pt);
  \fill[red] (2,2) circle (3pt);
  \fill[red] (4,1) circle (3pt);

  \draw[black, thick] (2,3) circle (5pt);
  \draw[black, thick] (4,2) circle (5pt);
  \draw[black, thick] (5,1) circle (5pt);

  \opendiamond{4}{2}
  \opendiamond{3}{4}
  \opendiamond{1}{5}

    \opensquare[.2]{1}{6}
    \opensquare[.2]{3}{5}
    \opensquare[.2]{4}{4}
    \opensquare[.2]{5}{2}
\end{tikzpicture}\end{center}
\caption{For the upset $U\subseteq \N^3$ of \cref{Ex:WzXz}, we show $M_{U_0}$ (blue) and $M_{U_1}\setminus M_{U_0}$ (red), as well as the $(x,y)$-projections of $X^1$ (diamonds), $W^0$ (squares), and $W^1$ (circles).}
\label{fig:WzXz}
\end{figure}
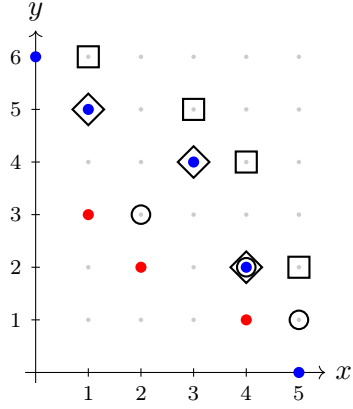
 \cref{Prop:Upset_to_Interval}\,(ii) yields an initial scaffold $P\subseteq U$ with \[\Ob P=M_U \cup W^0\cup W^1\cup X^1\] and the following relations.

Horizontal $X^1$-relations:
\begin{align*}
(1,3,1)&\leq (1,5,1),\\
(l,1)&\leq (3,4,1), \textup{ where $l$ is one of either $(1,3)$ or $(2,2)$}.
\end{align*}

Vertical $X^1$-relations:
\[(1,5,0)\leq (1,5,1),\quad (3,4,0)\leq (3,4,1),\quad (4,2,0)\leq (4,2,1).\]

 $W^0$-relations:
\begin{align*}
(0,6,0)&\leq (1,6,0)\geq (1,5,0)\leq (3,5,0) \geq (3,4,0)\\ (3,4,0)&\leq (4,4,0)\geq (4,2,0)\leq (5,2,0)\geq (5,0,0).
\end{align*}

$W^1$-relations:
\[(1,3,1)\leq(2,3,1)\geq  (2,2,1) \leq(4,2,1)\geq  (4,1,1)\leq(5,1,1)\geq  (5,0,0).\]

\end{example}

Our proof of \cref{Prop:Upset_to_Interval} depends on the next two lemmas. In what follows, all downsets are taken in $U$.

\begin{lemma}\label{Lem:Local__Connectivity}
For $p\in \N^3$, we have  $p\in E_U$ if and only if for some $i,j\in  \{1,2,3\}$ with $i\ne j$, we have
\begin{enumerate}
\item $\{p-\e_i,p-\e_j\}\subseteq \Down{p}$,
\item  $p-\e_i-\e_j \not\in \Down{p}$,
\item   $\{p-\e_i-\e_h,p-\e_j-\e_h\}\not\subseteq \Down{p}$, where $\{h\}=\{1,2,3\}\setminus \{i,j\}$.
\end{enumerate}
\end{lemma}

\begin{proof}
By definition, $p\in E_U$ if and only if $\Down{p}$ has at least two components.  Since $U$ is an upset, each component of $\Down{p}$ contains some point $p-\e_i$.
For $T$ as defined in \cref{Lem:Connectivity_of_G_and_Down_Z} and $i,j\in  \{1,2,3\}$ with $i\ne j$, we have that $p-\e_i$ and $p-\e_j$ belong to different components of $T$ if and only if conditions (1)-(3) are satisfied.  The result now follows from \cref{Lem:Connectivity_of_G_and_Down_Z}.
\end{proof}

\begin{lemma}\label{Lem:Local_Characterization_of_Wz_Xz}
For $p\in \N^3$ with $p_3=z$, we have
\begin{itemize}
\item[(i)] $p\in X^z$ if and only if
\begin{enumerate}
\item $p-\e_3\in \Down{p}$,
\item $p-\e_i\in \Down{p}$ for some $i\in \{1,2\}$, and
\item $p-\e_i-\e_3\not\in \Down{p}$ for each $i\in\{1,2\}$.
\end{enumerate}
\item[(ii)] $p\in W^z$ if and only if
\begin{enumerate}
\item $\{p-\e_1,p-\e_2\}\subseteq \Down{p}$,
\item $p-\e_1-\e_2\not\in \Down{p}$, and
\item $\{p-\e_1-\e_3,p-\e_2-\e_3\}\not\subseteq \Down{p}$.
\end{enumerate}
\end{itemize}
\end{lemma}

\begin{proof}
To prove (i), note that for $p=(m,z) \in \N^3$, $m\in M_{U_{z-1}}$ if and only if properties (1) and (3) hold.  Moreover, if $m\in M_{U_{z-1}}$, then $(2)$ holds if and only if $m\not\in M_{U_z}$.  It follows that $p\in X^z$ if and only if all three properties hold.

To prove (ii), consider the point $p=(m^i\vee m^{i+1},z)\in W^z$.  Note that $(m^i,z),(m^{i+1},z)\in \Down{p}$.  Since $(m^i,z)\leq p-\e_1<p$ and $(m^{i+1},z)\leq p-\e_2<p$, (1) holds.
  In addition, (2) holds because
  otherwise $m^i$ and $m^{i+1}$ would not be consecutively ordered by $x$-coordinate in $M_{U_z}$, a contradiction; indeed, if (2) does not hold, then there exists some element $m'\in M_{U_z}$ with $(m',z)\leq p-\e_1-\e_2$, and we have $m^i_1<m'_1<m^{i+1}_1$.  Note that since $U$ is an upset, (2) implies that $p-\e_1-\e_2-\e_3\not \in \Down{p}$.  We show by contradiction that (3) holds:
 If $\{p-\e_1-\e_3,p-\e_2-\e_3\}\subseteq \Down{p}$, then since also $p-\e_1-\e_2-\e_3\not\in \Down{p}$, we have $m^i\vee m^{i+1}\in V^{z-1}$, which contradicts the assumption that $p\in W^z$.

 Conversely, assume that $p=(m,z)\in \N^3$ satisfies properties (1)-(3) of part (ii) of the lemma statement.  Let
 \begin{align*}
 l^1=\min \{(a,m_2) \in U_z\mid a<m_1\},\\
 l^2=\min \{(m_1,b) \in U_z\mid b<m_2\}.
\end{align*}
By (1), the sets above are non-empty, so $l^1$ and $l^2$ are well defined.  By (2), $l^1$ and $l^2$ are consecutive elements of $M_{U_z}$, i.e., $l^1=m^i$ and $l^2=m^{i+1}$ for some $i$.  Moreover, $m=m^{i}\vee m^{i+1}$, so $m\in V^z$.  Property (3) implies that $m\not \in V^{z-1}$.  Hence $p\in W^z$.  This proves (ii).
\end{proof}

\begin{proof}[Proof of \cref{Prop:Upset_to_Interval}]
To prove (i), it suffices to show that for $p\in \N^3$ with $p_3=z$, we have $p\in W^z\cup X^z$ if and only if the three conditions of \cref{Lem:Local__Connectivity} hold.  By \cref{Lem:Local_Characterization_of_Wz_Xz}\,(ii),  $p\in W^z$ if and only if the conditions of \cref{Lem:Local__Connectivity} hold with $\{i,j\}=\{1,2\}$.  By \cref{Lem:Local_Characterization_of_Wz_Xz}\,(i), if $p\in X^z$, then the conditions of \cref{Lem:Local__Connectivity} hold with $\{i,j\}=\{1,3\}$ or $\{i,j\}=\{2,3\}$.  Conversely, suppose that the conditions of \cref{Lem:Local__Connectivity} hold for $\{i,j\}=\{1,3\}$.  We show that then $p\in W^z\cup X^z$: If $p-\e_2-\e_3\in \Down{p}$, then the conditions of \cref{Lem:Local__Connectivity} also hold for $\{i,j\}=\{1,2\}$, implying that $p\in W^z$.  If $p-\e_2-\e_3\not \in \Down{p}$, then the conditions of \cref{Lem:Local_Characterization_of_Wz_Xz}\,(i) hold, so $p\in X^z$.  A symmetric argument shows that if the conditions of \cref{Lem:Local__Connectivity} hold for $\{i,j\}=\{2,3\}$, then $p\in W^z\cup X^z$.
This completes the proof of (i).

To prove (ii), suppose $p=(m,z)\in X^z$.  Then \cref{Lem:Local_Characterization_of_Wz_Xz}\,(i) implies that $\Down{p}$ has a component \[A\subseteq\{(m,c)\in \N^3\mid c<z\}.\]
It follows that $P$ contains the relation $m_{\downarrow}<p$.  If in addition $p\not \in W^z$, then we claim that $\Down{p}$ has exactly one other component besides $A$, and that all points of this component have third coordinate $z$.
To see this, first note that by \cref{Lem:Local_Characterization_of_Wz_Xz}\,(i), $\Down{p}$ has at least one other component besides $A$, which must contain either $p-\e_1$ or $p-\e_2$.  Moreover, since there is no $l<m$ with $(l,z-1)\in \Down{p}$, each such component contains only points with third coordinate $z$.  To verify the claim, it thus suffices to check that  $\Down{p}$ has a single component containing either $p-\e_1$ or $p-\e_2$.  We consider three cases:
\begin{itemize}
\item If $\{p-\e_1,p-\e_2\}\subseteq \Down{p}$, then since $p\in X^z\setminus W^z$, \cref{Lem:Local_Characterization_of_Wz_Xz} implies that $p-\e_1-\e_2\in \Down{p}$.
\item If $p-\e_1\in \Down{p}$ and $p-\e_2\not \in \Down{p}$, then every element of $\Down{p}$ with third coordinate $z$ is contained in the set $\{(a,m_2,z)\in \N^3\mid a<m_1\}$.
\item Symmetrically, if $p-\e_2\in \Down{p}$ and $p-\e_1\not \in \Down{p}$, then every element of $\Down{p}$ with third coordinate $z$ is contained in the set  $\{(m_1,b,z)\in \N^3\mid b<m_2\}$.
\end{itemize}
In each case, $\Down{p}$ has a single component containing either $p-\e_1$ or $p-\e_2$, so the claim holds.
For $(l,z)$ an arbitrary point of $M_{U}$ with $l<m$, the claim implies that $(l,z)<p$ is a valid choice of relation in $P$.  Since $\Down{p}$ has two components, no other relations with upper point $p$ need to be added to $P$.

Next, suppose that $p=(m^i\vee m^{i+1},z)\in W^z$.  Since  $(m^i,z)\leq p-\e_1$ and $(m^{i+1},z)\leq p-\e_2$, the proof of (i) shows that these points lie in separate components of $\Down{p}$.  Hence $m^i_{\downarrow}<p$ and $m^{i+1}_{\downarrow}<p$ are valid choices of relations in $P$.  If in addition $p\not\in X^z$, then by \cref{Lem:Local_Characterization_of_Wz_Xz}\,(i) one of the following holds:
\begin{itemize}
\item $p-\e_3\not\in U$,
\item $p-\e_2-\e_3\in U$, hence $p-\e_2$ and $p-\e_3$ lie in the same component of $\Down{p}$,
\item $p-\e_1-\e_3\in U$, hence $p-\e_1$ and $p-\e_3$ lie in the same component of $\Down{p}$.
\end{itemize}
In each case, \cref{Lem:Connectivity_of_G_and_Down_Z} implies that $\Down{p}$ has exactly two components, so no other relations with upper point $p$ need to be added to $P$.

Finally, if $p=(m,z)\in X^z\cap W^z$ with $m=m^i\vee m^{i+1}$, then the case $p\in W^z$ treated above yields two relations $m^{i}_{\downarrow}<p$ and $m^{i+1}_{\downarrow}<p$, while the case $p\in X^z$ yields a third relation $m_{\downarrow}<p$.  As $p\in X^z$, there is no $l<m$ with $(l,z-1)\in \Down{p}$, so $m^i_{\downarrow}$ and $m^{i+1}_{\downarrow}$ both have third coordinate $z$.  Thus, $m^i_{\downarrow}$, $m^{i+1}_{\downarrow}$, and $m_{\downarrow}$ lie in different components of $\Down{p}$, so all three relations can be included in $P$.  Since each component of $\Down{p}$ must contain some point $p-\e_i$, it is not possible for $\Down{p}$ to have more than three components.  Therefore, no other relations with upper point $p$ need to be added to $P$.  This completes the proof of (ii).
\end{proof}

Our algorithm will compute an initial scaffold $P$ of $U$ as in \cref{Prop:Upset_to_Interval}.  The algorithm will identify points of $W^z$ using the following characterization of $W^z$ in terms of $X^z$:

\begin{lemma}\label{Lem:W^z_characterization_in_terms_of_X^z}
For $m^i\vee m^{i+1}\in V^z$, we have $(m^i\vee m^{i+1},z) \in W^z$ if and only if either
\begin{enumerate}
\item  $(m^{i+1},z) \in M_{U}$ and there is no $(m^{i+1}_1,b,z)\in X^z$ with $b<m^{i}_2$, or
\item $(m^{i},z) \in M_{U}$ and there is no $(a,m^i_2,z)\in X^z$ with  $a<m^{i+1}_1$.
\end{enumerate}
\end{lemma}

\begin{proof}
Note that $(m^i\vee m^{i+1},z) \in W^z$ if and only if $(m^i\vee m^{i+1}) \not\in V^{z-1}$.  Moreover, $(m^i\vee m^{i+1})\in V^{z-1}$ if and only if both of the following are true:
\begin{enumerate}
\item  $m^{i+1} \in M_{U_{z-1}}$ or there exists $(m^{i+1}_1,b,z)\in X^z$ with $b<m^{i}_2$,
\item  $m^{i} \in M_{U_{z-1}}$ or there exists $(a,m^i_2,z)\in X^z$ with $a<m^{i+1}_1$.
\end{enumerate}
Moreover, for $j\in \{i,i+1\}$, we have $(m^{j},z) \in M_{U}$ if and only if $m^{j} \not \in M_{U_{z-1}}$.  Thus, by contraposition and de Morgan's laws, we have $(m^i\vee m^{i+1}) \not\in V^{z-1}$ if and only if one of the two conditions in the lemma statement holds.  The result follows.
\end{proof}

\begin{example}
In \cref{Ex:WzXz}, we have seen that $(1,6,1)\not\in W^1$.  To illustrate \cref{Lem:W^z_characterization_in_terms_of_X^z}, we check that $(1,6,1)$ satisfies neither condition (1) nor (2) of the lemma.  First, letting $m^1=(0,6)$ and $m^2=(1,3)$, note that $m^1,m^2\in M_{U_1}$ and $(1,6)=m^1\vee m^2$.  Condition (1) fails because the point $(1,5,1)\in X^1$ is of the form $(m^2_1,b,1)$, where $b=5<6=m^1_2$.  Condition (2) fails because $(m^1,1)=(0,6,1) \not\in M_{U}$.
\end{example}

The following technical observation will be helpful for understanding how our algorithm computes relations of $P$ whose upper points lie in $X^z\cap W^z$.

\begin{remark}\label{Rem:relations_of_Xz_cap_Wz}
If $p=(m^i\vee m^{i+1},z)\in X^z\cap W^z$, then for any point $(l,z)\in M_U\cap \Down{p}$, we must have $l\in \{m^i,m^{i+1}\}$, since $p-\e_1-\e_2\not\in \Down{p}$ by \cref{Lem:Local_Characterization_of_Wz_Xz}\,(ii).  Hence, the relation $(l,z)<p$ is a $W^z$-relation.
\end{remark}

\begin{example}
In \cref{Ex:WzXz}, we have seen that $p=(4,2,1)\in X^1\cap W^1$.  To illustrate \cref{Rem:relations_of_Xz_cap_Wz}, note that the set of points in $M_U\cap \Down{p}$ with third coordinate 1 is $\{(2,2,1),(4,1,1)\}$.  We have seen in \cref{Ex:WzXz} that both $(2,2,1)<p$ and $(4,1,1)<p$ are $W^1$-relations.
\end{example}

To prepare for our algorithm, we recall that a self-balancing binary search tree $\mathcal T$, such as a red-black tree \cite[Chapter 13]{cormen2022introduction}, stores an ordered set in a way that supports binary searches, insertions, and deletions, each in $O(\log |\mathcal T|)$ time, where $|\mathcal T|$ is the number of elements stored.   In addition,  $\mathcal T$ supports \emph{range queries}: Given keys $r$ and $s$, it returns the set $[r,s]\coloneqq\{t\in \mathcal T\mid r\leq t\leq s\}$ as an ordered list in time $O(\log |\mathcal T|+|[r,s]|)$; see \cite[Exercise 12.2-8]{cormen2022introduction}.

\begin{proof}[Proof of \cref{Thm:Computing_Hull_Skeleton_Intervals}\,(i)]
Given $M_U=M_Q$, our algorithm first computes an initial scaffold $P$ of $U$ as in \cref{Prop:Upset_to_Interval} by sweeping through the planar slices $U_{z}$ of $U$.  We represent $P$ in memory as a directed graph, where vertices are the points of $P$, edges are the relations, and each vertex maintains pointers to its incident edges.
 Initially, we set $P\leftarrow M_U$, with no non-identity relations.  As the algorithm proceeds, it adds points and relations to $P$.

Let \[T^z=\{m_{\downarrow}\mid m\in M_{U_{z}}\}.\]
Note that $M_{U_z}$ is the image of $T^z$ under coordinate projection onto the $xy$-plane.
We maintain $T^z$, ordered by $x$-coordinate, using a self-balancing binary search tree  $\mathcal T$.  Initially, we set $\mathcal{T}\leftarrow \emptyset$.  As we modify $\mathcal T$ during the course of our algorithm, $\mathcal T$ always stores a set of points whose image under coordinate projection onto the $xy$-plane is pairwise incomparable.  Thus, the order of $\mathcal T$ by $x$-coordinate is the same as the opposite of the order by $y$-coordinate.

Let
\begin{align*}
Y^z&=\{m\in \N^2\mid (m,z)\in M_U\},\\
   Z&=\{z\in \N\mid Y^z\ne \emptyset\}.
 \end{align*}
Note that if $z\not \in Z$, then $X^z=W^z=\emptyset$.  Thus, by \cref{Prop:Upset_to_Interval}\,(i), to compute $E_U$, it suffices to compute $X^z$ and $W^z$ for each $z\in Z$.  Also  note that for consecutive elements $z'<z$ of $Z$, we have $M_{U_{z-1}}=M_{U_{z'}}$, hence $T^{z-1}=T^{z'}$, and that if $z<\min(Z)$, then $T^{z}=\emptyset$.

We iterate through the elements of $Z$ in increasing order.    Inductively, for $z\in Z$, assume that $\mathcal T=T^{z-1}$ has been computed.  We store $Y^z$ in a list ordered by $x$-coordinate.  Using $T^{z-1}$ and $Y^z$, we compute $X^z$, $T^z$, $W^z$, the $X^z$-relations, and the $W^z$-relations.  This computation consists of two phases.  In the first phase, we compute $X^z$, $T^z$, and a set of relations of $P$ which includes all of the $X^z$-relations and some of the $W^z$-relations.  In the second phase, we compute $W^z$ and the remaining $W^z$-relations.

\subsubsection*{First Phase} The first phase proceeds as follows:
We initialize an empty self-balancing binary search tree that will store the points of $X^z$, ordered by $x$-coordinate; as we discover points of $X^z$, we insert them into this tree.  We iterate through the points of $Y^z$ in any order.  For each $m\in Y^z$, we compute the set
\[S^m\coloneqq \{m'_{\downarrow}\in \mathcal{T}\mid m<m'\};\]
the details of this computation are explained below.  We then insert each point of \[\{(m',z) \mid m'_{\downarrow} \in S^m\}\] into both $X^z$ and $P$, and also add both of the relations $m'_{\downarrow}<(m',z)$ and $(m,z)<(m',z)$ to $P$.
(Note that the former relation is a vertical $X^z$-relation, while the latter is a  horizontal $X^z$-relation if $(m',z)\not \in W^z$ and, by \cref{Rem:relations_of_Xz_cap_Wz}, a $W^z$-relation if $(m',z)\in W^z$.)
Using insertion and deletion operations, we then set
\[\mathcal{T}\leftarrow (\mathcal{T}\setminus S^m) \cup \{(m,z)\}.\]

To explain how the set $S^m$ is computed, let $L$ be the prefix of $\mathcal{T}$ consisting of points with
$y$-coordinate at least $m_2$, and let $R$ be the suffix of $\mathcal{T}$ consisting of points with $x$-coordinate at least $m_1$. Note that $S^m=L\cap R$, as no point of $\mathcal{T}$ projects to $m$ prior to the computation of $S^m$.    We determine the minimum point $R_{\min}$ of $R$ (null if $R$ is empty) and the maximum point $L_{\max}$ of $L$ (null if $L$ is empty) by a pair of binary searches on $\mathcal{T}$.  We then obtain $S^m=L\cap R=[R_{\min},L_{\max}]$ by a range search on $\mathcal T$.

After iterating over all points of $Y^z$, we have $\mathcal{T}=T^z$.  In addition, we have fully computed $X^z$ and the $X^z$-relations.  This completes the first phase.

\subsubsection*{Second Phase} In the second phase, we compute $W^z$ and the remaining $W^z$-relations by applying \cref{Lem:W^z_characterization_in_terms_of_X^z}.  As above, we write \[M_{U_z}=(m^1,\ldots,m^k).\]
\cref{Lem:W^z_characterization_in_terms_of_X^z} tells us in particular that for each $(m^i\vee m^{i+1},z)\in W^z$, either $m^i\in Y^z$ or $m^{i+1}\in Y^z$.  We can therefore discover all elements of $W^z$ by iterating through all $m^j\in Y^z$ and checking whether each of $m^{j-1}\vee m^{j}$ and $m^{j}\vee m^{j+1}$ is an element of  $W^z$.  We call these two checks \emph{leftward search} and \emph{rightward search}, respectively.  We perform these checks using the characterization of $W^z$ given in \cref{Lem:W^z_characterization_in_terms_of_X^z}.

In more detail,
we proceed as follows: We initialize an integer variable $\textsc{LastRightwardDiscovery}=0$, which we use to avoid processing an element $(m^{i}\vee m^{i+1},z)\in W^z$ twice when  $m^i,m^{i+1}\in Y^z$.  For each $m^j\in Y^z$ in increasing order, we do the following:
\begin{itemize}
\item{} \textbf{Leftward search}: If $j>1$ and $\textsc{LastRightwardDiscovery}\ne j-1$,
then we check via binary search whether there exists a point $(m^j_1,b, z)\in X^z$ with  $b<m^{j-1}_2$.
If not, then \[p\coloneqq (m^{j-1}\vee m^j,z)\in W^z,\] so we do the following:
\begin{itemize}
\item If $p\not\in X^z$, then we add $p$ to $P$, along with the relations $m_{\downarrow}^{j-1}<p$ and $m^{j}_{\downarrow}<p$.  If $p\in X^z$, then one of these two relations was already added to $P$ in the first phase of the algorithm; we add the other to $P$.
\end{itemize}
\item{} \textbf{Rightward search}: If $j<k$, then we check whether there exists a point $(a,m^{j}_2,z)\in X^z$ with $a<m^{j+1}_1$.
If not, then \[p\coloneqq (m^{j}\vee m^{j+1},z)\in W^z,\] so we do the following:
\begin{itemize}
\item If $p\not\in X^z$, then we add $p$ to $P$, along with the relations $m_{\downarrow}^{j}<p$ and $m^{j+1}_{\downarrow}<p$.  If $p\in X^z$, then one of these two relations was already added to $P$ in the first phase of the algorithm; we add the other to $P$.
\item We set $\textsc{LastRightwardDiscovery}\leftarrow j$.
\end{itemize}
 \end{itemize}
We claim that after iterating over all points of $Y^z$, each element $p=(m^i\vee m^{i+1},z)$ of $W^z$ is discovered.  To see this, note that if $p$ satisfies condition (2) of \cref{Lem:W^z_characterization_in_terms_of_X^z}, then the rightward search at index $j=i$ discovers that $p\in W^z$ and sets $\textsc{LastRightwardDiscovery}\leftarrow i$.  If $p$ satisfies condition (1) but not condition (2) of \cref{Lem:W^z_characterization_in_terms_of_X^z}, then $\textsc{LastRightwardDiscovery}$ is never set to $i$, so leftward search at index $j=i+1$ discovers that $p\in W^z$.

This completes the computation of $P$.

\subsubsection*{Computation of $P\cap Q$}
To explain how we compute $P\cap Q$ from $P$, first note that since $Q$ is an interval, it is a full subposet of $U$.  Thus, since $M_Q=M_U$, to compute $P\cap Q$ from $P$, it suffices to compute $\Ob(P\cap Q)=M_U\cup (E_U\cap Q)$, recording all relations of $P$ with upper point in $E_U\cap Q$.

To compute $E_U\cap Q$ from $P$, assume first that $Q$ is specified as input to our algorithm by its upset presentation $(M_Q,M_{Q'})$.
We sweep through the planar slices
\[{Q'_z}\coloneqq\{(x,y)\in \N^{2} \mid (x,y,z)\in Q'\},\] computing each $M_{Q'_{z}}$ from $M_{Q'_{z-1}}$ exactly as we computed $T^z$ from $T^{z-1}$ above.  (Here, we needn't store the third coordinates of elements of $\mathcal{T}$, which were needed only to compute the relations of $P$, so $\mathcal T$ stores the sets $M_{Q'_{z}}$ directly.)  For $p=(x,y,z) \in E_U$, note that
$p\not \in Q$ if and only if $p\in Q'$, and that $p\in Q'$ if and only if $(x,y)\geq m$ for some $m\in M_{Q'_{z}}$.  The last condition can be checked using at most two binary searches on $M_{Q'_z}$.

If $Q$ is finite and instead specified by its extrema, a dual variant of the above procedure computes $E_U\cap Q$ from $P$:  For any poset $A$, let $N_A$ denote the set of maxima of $A$.  Let $D=\DownCL{N_Q,\N^3}$, and for $z\in \N$, let \[D_{z}\coloneqq\{(x,y)\in \N^{2} \mid (x,y,z)\in D\}.\]  We sweep through the planar slices $D_{z}$ in decreasing order of $z$, computing each $N_{D_{z}}$ from $N_{D_{z+1}}$.   For $p=(x,y,z) \in E_U$, note that $p\in Q$ if and only if $(x,y)\leq w$ for some $w\in N_{D_{z}}$.  As above, the latter condition can be checked by at most two binary searches on $N_{D_z}$.  This completes the specification of our algorithm.

\subsubsection*{Complexity}
We now turn to the complexity analysis.  We first consider the cost of computing $P$.  Note that \[M_Q=M_U=\bigsqcup_{z} \{(m,z)\mid m\in Y^z\},\] so $\sum_{z} |Y^z|=|M_Q|=O(n)$.
For each $m\in Y^z$, we perform $O(1)$ binary searches on $\mathcal{T}$ and $O(1)$ binary searches on $X^z$.  Since $|\mathcal{T}|\leq |M_Q|$ at all times and $|X^z|\leq |M_Q|$ for all $z$, the total time required by these binary searches across all levels $z$ is $O(|M_Q|\log |M_Q|)=O(n\log n)$.  In addition, across all levels $z$, each element of $M_Q$ is inserted into $\mathcal{T}$ exactly once and removed from $\mathcal{T}$ at most once.  Thus, the cost of all such insertions and deletions is again $O(n\log n)$.   By \cref{Prop:Upset_to_Interval}\,(i), each element inserted into some set $X^z$ belongs to $E_U$, and each element of $E_U$ is inserted into a set $X^z$ at most once, across all levels $z$.   Since $|E_U|=O(|M_Q|)$ by \cref{Thm:Size_of_Initial}, the cost of all such insertions is once again $O(n\log n)$.   For each $m\in Y^z$, the associated range search on $\mathcal T$ takes time $O(\log |\mathcal T|+s_m)=O(\log |M_Q|+s_m)$, where $s_m$ is the number of elements removed from $\mathcal T$ when we process $m$.  Since \[\sum_{z}\sum_{m\in Y^z} s_m \leq |M_Q|,\] the total cost of all range searches is $O(n\log n)$.  Besides these operations, computing $P$ requires $O(1)$ time per element of $Y^z$.  Hence the total time to compute $P$ is $O(n\log n)$.

To compute $\Ob(P\cap Q)$, we must compute $M_{Q'_z}$ (or $N_{D_z}$) at all values of $z$ for which there is an element of $M_{Q'}$ (or $N_{Q}$) with third coordinate $z$; this takes time $O(n\log n)$.  Besides this, we must perform $O(|E_U|)$ binary searches on lists of size $O(n)$.  By \cref{Thm:Size_of_Initial}, we have $|E_U|=O(|M_Q|)=O(n),$ so the total time required by all of these binary searches is $O(n\log n)$.
As each element $p\in E_U$ maintains pointers to the relations $m<p$ in $P$, the cost of computing $\hom(P\cap Q)$ is dominated by the cost of computing $\Ob(P\cap Q)$.  Thus, the total time to compute $P\cap Q$ from $P$ is $O(n\log n)$.
\end{proof}

\begin{remark}
In view of \cref{Essential_and_Betti_1}\,(ii), the algorithm underlying our proof of \cref{Thm:Computing_Hull_Skeleton_Intervals}\,(i) specializes to an $O(n \log n)$-time algorithm for computing the support of the first Betti numbers of a monomial ideal $J\subseteq R^3$, where $n$ is the size of a minimal generating set of $J$.  A slight modification of the algorithm yields the first Betti numbers themselves, using the fact that for $J=\hat U$ and $p\in \N^3$ with $p_3=z$,
\[\beta^J_1(p)=
\begin{cases}
2 &\textup{if }p\in X^z\cap W^z,\\
1 &\textup{if }p\in (X^z\setminus W^z)\cup (W^z\setminus X^z),\\
0 &\textup{otherwise},
\end{cases}\]
which follows immediately from \cref{Essential_and_Betti_1}\,(i) and \cref{Prop:Upset_to_Interval}.  We expect that the slice-wise approach of our algorithm also extends to compute a minimal resolution of $J$ (and in particular, the second Betti numbers of $J$), in the same asymptotic time.  An outline of a different algorithm for this has previously been given in \cite[Section 3]{miller2005combinatorial}, without explicit complexity bounds.
\end{remark}

\subsection{Algorithm for Intervals in $\N^d$}\label{Sec:Compute_Hull_Skeleta_Nd}
We next prove \cref{Thm:Computing_Hull_Skeleton_Intervals}\,(ii) by giving an algorithm to compute an initial scaffold of an interval in $\N^d$, for arbitrary $d$.  As in the $d=3$ case, our strategy is to first compute an initial scaffold $P$ of $U=\Up{Q,\N^d}$ and then compute $P\cap Q$, which is an initial scaffold  of $Q$ by \cref{Prop:Hull_Skeleton_as_Subset}.

Let $T\subseteq U$ be the full subposet consisting of points of $M_U$ and their pairwise joins.  To compute $P$, we use the following structural result.

\begin{proposition}\label{Prop:Joins_Suffice}
A poset $P\subseteq U$ is an initial scaffold of $U$ if and only if it is an initial scaffold of $T$.
\end{proposition}

\begin{proof}
It is easily checked that any point of $E_U$ is a pairwise join of points in $M_U$.  (Equivalently, any point in the support of $\beta_1^{\hat U}$ is such a join, which one sees immediately by considering the Taylor resolution of $\hat U$ \cite{miller2005combinatorial}.)  Hence $E_U\subseteq T$.  We claim that for any $t\in T$, the inclusion $\Down{t,T}\hookrightarrow \Down{t,U}$ induces a bijection on components.  The claim implies that $E_U=E_T$, and hence that $I_U=I_T$.  Since $M_U=M_T$, it then follows that $U$ and $T$ have the same initial scaffolds, as desired.

It remains to prove the claim.   Our argument is quite similar to the proof of \cref{Lem:Connectivity_of_G_and_Down_Z}.  Since $U$ is well-founded, every component of $ \Down{t,U}$ contains a minimal element.  Since $M_T=M_U$,  the inclusion $\Down{t,T}\hookrightarrow \Down{t,U}$ therefore induces a surjection on components.  We need to show that the induced map is also an injection.  Suppose $w,y\in T$, with $w$ and $y$ contained in the same component of $\Down{t,U}$.  It suffices to show that $w$ and $y$ are contained in the same component of $\Down{t,T}$.   By assumption, there exists a path $w=q_1,\ldots,q_k=y$ from $w$ to $y$ in $\Down{t,U}$.  We may assume without loss of generality that for each  triple $q_{i-1},q_i,q_{i+1}$ of consecutive points in the path, either \[q_{i-1}< q_i> q_{i+1}\quad \textup{or} \quad q_{i-1}> q_i< q_{i+1}.\]  We call $q_i$ an upper point in the former case and a lower point in the latter case.  Moreover, we may assume without loss of generality that if $q_1$ is non-minimal and $k>2$, then $q_2$ is a lower point and, similarly, that if $q_k$ is non-minimal and $k>2$, then $q_{k-1}$ is a lower point.  Replacing each non-minimal lower point with a minimal point below it, and then replacing each upper point $q_{i}$ with the join $q_{i-1}\vee q_{i+1}$, we obtain a path from $w$ to $y$ in $\Down{t,T}$.  Thus, $w$ and $y$ lie in the same component of $\Down{t,T}$, as we wanted.
\end{proof}

\begin{proof}[Proof of {\cref{Thm:Computing_Hull_Skeleton_Intervals}}\,(ii)]
By \cref{Prop:Joins_Suffice}, to compute an initial scaffold $P$ of $U=\Up{Q}$, it suffices to compute an initial scaffold of $T$.  To do this, we initialize $P\leftarrow M_U$, with no non-identity relations.  We then compute $T\setminus M_T$,  ordered in a way compatible with the partial order induced by $\N^d$, e.g., lexicographically.

We next iterate through the points of $T\setminus M_T$ in order, adding each $t\in E_T$ to $P$, together with the required relations with upper point $t$.  To elaborate, let $P^t$ be the part of $P$ computed so far when we come to $t\in T\setminus M_T$.  By \cref{Lem:Hull_Skeleton_Path_Components}, to decide whether $t\in E_T$ and, if so, to compute the relations of $P$ with upper point $t$, it suffices to compute the components of ${P^t\cap \Down{t,T}}=P\cap \Down{t,T}$ rather than those of $\Down{t,T}$.  We compute $P^t\cap \Down{t,T}$ by comparing each point of $P^t$ to $t$.  Similarly to our algorithm for computing initial scaffolds of general finite posets given in the proof of \cref{Thm:Compute_Hull_Skeleton}, we then compute the components of $P^t\cap \Down{t,T}$  via depth-first search.
If $P^t\cap \Down{t,T}$ has more than one component, then we add $t$ to $P$, and for each component $A$ of $P^t\cap \Down{t,T}$, we choose a minimal element $m\in A$ and add the relation $m<t$ to $P$.

For each point $t\in T\setminus M_T$ considered in the above algorithm, identifying the components of $P^t\cap \Down{t,T}$ requires time $O(|\hom P|)=O(|M_Q|^2)=O(n^2)$, where the first equality follows from \cref{Thm:Size_of_Initial}.  At most $d$ relations are added to $P^t$ per point $t$, so adding these relations requires constant time.   Since $T$ has $O(|M_Q|^2)=O(n^2)$ elements,
the total time to compute the initial scaffold $P$ is $O(n^4)$.

 Given $P$, we can straightforwardly compute $\Ob(P\cap Q)$, and hence $P\cap Q$, in time $O(n^3)$.  To explain, if the input to our algorithm is the upset presentation $(M_Q,M_{Q'})$, then
 for each $p\in E_U$, checking whether $p\in Q$ amounts to checking whether $p\geq m$ for some $m\in M_{Q'}$, as in the proof of \cref{Thm:Computing_Hull_Skeleton_Intervals}\,(i).   Since $|M_{Q'}|=O(n)$ and $|E_U|=O(n^2)$ by \cref{Thm:Size_of_Initial}, the total time required for all such comparisons is $O(n\cdot n^2)=O(n^3)$.  If $Q$ is instead given by the extrema of $Q$, then the computation of $P\cap Q$ is similar.  Thus, the total time to compute an initial scaffold of $Q$ is $O(n^4)$.
\end{proof}

\section{Proofs of Results on Limit Computation}\label{Sec:Limit_Comp_Proofs}
In this section, we complete the proofs of our results on limit computation by proving \cref{Prop:Chain_Complex_Compute_Matrix_Rep} and
\cref{Thm:Chain_complex_limit}\,(i).

\subsection{Proof of \cref{Prop:Chain_Complex_Compute_Matrix_Rep}} \label{Sec:From_Chain_Complex_to_Matrix_Rep}

Recall the statement of the proposition:

\begin{repproposition}{Prop:Chain_Complex_Compute_Matrix_Rep}
Suppose we are given
\begin{enumerate}
\item a finite poset $P$, specified by its sets of elements and relations,
\item a poset $Q$, represented in such a way that pairs of elements can be compared in constant time,
\item a functor $j \colon P\to Q$, specified by its values on objects, and
\item a $(Q,r)$-complex with homology $H$.
\end{enumerate}
Then we can compute a matrix representation of $H \circ j$ in time $O(|\hom P|r^{\omega})$.
\end{repproposition}

We prepare for the proof by recalling some elementary linear algebra results, which we state without proof.

\begin{definition}
For $\kbb$-vector spaces $U\subseteq V$, we say a set of vectors $C$ in $V$ is a \emph{basis extension} of $U\subseteq V$ if $B\sqcup C$ is a basis of $V$ whenever $B$ is a basis  of $U$.
\end{definition}

\begin{lemma}
If $C$ is a basis extension of $U\subseteq V$, then the quotient map $V\to V/U$ sends $C$ to a basis $\bar C$ of $V/U$.
\end{lemma}

For a linear map $\alpha\colon V\to V'$ and ordered bases $B$ and $B'$ for $V$ and $V'$, we write the matrix representation of $\alpha$ with respect to these bases either as $[\alpha]^{B',B}$, or more simply as $[\alpha]$.

\begin{lemma}\label{lem:matrix_reps_of_induced_maps}
  Given a linear map $\alpha\colon V\to V'$ with $\alpha(U)\subseteq U'$, let $\bar \alpha\colon V/U\to V'/U'$ denote the induced map.  Let $B,B'$ be ordered bases for $U,U'$, and let
  $C,C'$ be ordered basis extensions of $U\subseteq V$, $U'\subseteq V'$.  Then the matrix $[\bar \alpha]^{\bar C',\bar C}$ is the block of the matrix $[\alpha]^{B'\sqcup C',B\sqcup C}$ consisting of rows and columns indexed by $C'$ and $C$, respectively.
 \end{lemma}

\begin{definition}\mbox{}
\begin{itemize}
\item[(i)]
Given a non-zero vector $c=(c_1,\ldots,c_l)\in \kbb^l$, the \emph{pivot index} of $c$ is the smallest index $j$ such that $c_j\ne 0$.
\item[(ii)] Given a finite-dimensional vector space $V$ with ordered basis $A$ and $v\in V$, let $[v]\in \kbb^{|A|}$ be the representation of $v$ with respect to $A$.
We say that an ordered set of non-zero vectors  $C=(v_1,\ldots,v_l )$ in $V$ is in \emph{echelon form} (with respect to $A$) if the pivot indices of the vectors $[v_k]$ are strictly increasing with $k$.
\end{itemize}
\end{definition}

\begin{proof}[Proof of \cref{Prop:Chain_Complex_Compute_Matrix_Rep}]
Recalling \cref{Def:QrComplex}, suppose that the given matrices $[f]$ and $[g]$ represent $f\colon X\to Y$ and $g\colon Y\to Z$ with respect to ordered bases $A^X,A^Y,A^Z$ for $X,Y,Z$.  For each $p\in Q$, these bases induce ordered bases $A^X_p,A^Y_p,A^Z_p$ for $X_p,Y_p,Z_p$.   The matrix representation $[f_p]$ of $f_p$ with respect to the induced bases is obtained from $[f]$ by taking the submatrix consisting of all rows and columns with labels $\leq p$.  Since we assume that pairs of elements of $Q$ can be compared in constant time, deciding whether a column or row of $[f]$ belongs to $[f_p]$ also takes constant time.  The same is true for $[g_p]$.

Computing a matrix representation of $H\circ j$ amounts to computing
\begin{enumerate}
\item a basis of $H_p$ for each $p\in j(\Ob P)$ and
\item the corresponding matrix representation of the map $H_{pq}$, for each  $p\leq q\in j(\hom P)$.
\end{enumerate}

We do this as follows:  First, for each $p\in j(\Ob P)$, we use Gaussian elimination with backsolve to compute a basis $B_p$ for $\im f_p$, a basis extension $C_p$ for $\im f_p \subseteq \ker g_p$, and a basis extension $D_p$ for $\ker g_p\subseteq Y_p$, where all elements are given as linear combinations of the basis $A^Y_p$.  Our algorithm selects the basis $\bar C_p$ for $H_p$.

We now explain this computation in detail.  To begin, we compute an ordered basis $A^{\ker g}_p$ of $\ker g_p$ in echelon form (with respect to $A^Y_p$).  We do this by applying Gaussian elimination with backsolve to $[g_p]$.   Writing $A^Y_p=\{a_1,\ldots,a_l\}$, we take
\[D_p=\{a_j\in A^{Y}_p \mid j \textup{ is not a pivot index of $[b]$ for any } b\in A^{\ker g}_p \}.\]
Multiplying $[f_p]=[f_p]^{A^Y_p,A^X_p}$ by a change-of-basis matrix on the left, we obtain the matrix $[f_p]^{(A^{\ker g}_p\sqcup D_p),A^X_p}$.  Let $\gamma_p\colon X_p\to \ker g_p$ denote the co-restriction of $f_p$ to $\ker g_p\subseteq Y_p$, i.e., $\gamma_p(v)=f_p(v)$ for all $v\in X_p$.  We obtain the matrix $[\gamma_p]^{A^{\ker g}_p,A^X_p}$ as the top block of $[f_p]$.  Performing Gaussian elimination on the transpose of this matrix yields a basis $B_p$ for $\im \gamma_p$ in echelon form with respect to $A^{\ker g}_p $.  We now compute the basis extension $C_p$ for $\im f_p\subseteq \ker g_p$ in the same way we computed $D_p$.  By the way $C_p$ is constructed, we have $C_p\subseteq A^{\ker g}_p$, so we immediately have expressions for the elements of $C_p$ in the basis $A^Y_p$.  To express the elements of the basis $B_p$ in the basis $A^Y_p$, it suffices to perform a single matrix multiplication.

For each relation $p\leq q$ in $j(\hom P)$, we need to compute the matrix $[H_{pq}]^{\bar C_q,\bar C_p}$.  To do this, first we form the matrix $[Y_{pq}]^{A^Y_{q},A^Y_{p}}$.  Note that this is the binary matrix whose entry at index $(y,z)$ is 1 if and only if the $y^{\mathrm{th}}$ element of $A^Y_{q}$ and the $z^{\mathrm{th}}$ element of $A^Y_{p}$ are induced by the same element of $A^Y$.  We then multiply this matrix by change-of-basis matrices on both the left and right to obtain the matrix \[W\coloneqq  [Y_{pq}]^{(B_q\sqcup C_q\sqcup D_q),(B_p\sqcup C_p\sqcup D_p)}.\]  It follows from \cref{lem:matrix_reps_of_induced_maps} that $[H_{pq}]^{\bar C_q,\bar C_p}$ is the block of $W$ with columns indexed by $C_p$ and rows indexed by $C_q$.

It remains to show that this algorithm runs in time $O(|\hom P|r^{\omega})$.  For each element or relation of $\im j $, the algorithm performs Gaussian elimination with backsolve a constant number of times and also performs a constant number of matrix multiplications.  Each matrix considered has dimension $O(r)\times O(r)$, so each instance of Gaussian elimination or matrix multiplication requires time $O(r^{\omega})$; see \cref{Omega_Linear_Solve}.  Since such operations dominate the cost of the algorithm, it follows that our algorithm runs in time $O(|\hom P|r^{\omega})$.
\end{proof}

\subsection{Proof of \cref{Thm:Chain_complex_limit}\,(i) }\label{Sec:Computing_Limit_Z^2}
Let us recall the statement of \cref{Thm:Chain_complex_limit}\,(i):

\begin{reponecorollary}{Thm:Chain_complex_limit}
Given an $(\N^2,r)$-complex with homology $H$ and the upset presentation of an interval $i\colon Q\hookrightarrow \N^2$ with $\|Q\|=n$, we can compute $\lim  (H\circ i)$ in time $O(n\log n+r^3+nr^2)$.
\end{reponecorollary}

Our proof will use the following decomposition result  from classical quiver representation theory, due to Gabriel.
\begin{theorem}[\cite{gabriel72,carlsson2010zigzag}]\label{Thm:Zigzag_decomposition}
Given a  zigzag poset $P$ (see \cref{Def:Zigzag_Poset}) and a pointwise finite-dimensional $P$-module $\fun$, there is a unique finite multiset $\mathcal B(\fun)$ of intervals in $P$ such that
\begin{equation}\label{Eq:Zigzag_Decomp}
\fun\cong \bigoplus_{I\in \mathcal B(\fun)} \kbb^{I},
\end{equation}
where the $\kbb^{I}$ are the interval modules of \cref{Def:Int_Mod}.
\end{theorem}
We call an isomorphism as in \cref{Eq:Zigzag_Decomp} an \emph{interval decomposition of $\fun$}.

\begin{proposition}\label{Prop:Limit_of_Zigzag}
In the setting of \cref{Thm:Zigzag_decomposition}, we have
\begin{equation}\label{eq:limits_products_commute}
\lim \fun\cong \bigoplus_{I\in \mathcal B(\fun)} \lim \kbb^{I},
\end{equation}
where
\[
\lim \kbb^{I}=\begin{cases}
 \kbb & \textup{ if $I$ is a downset of $P$},\\
             0 &\textup{otherwise.}
             \end{cases}
\]
\end{proposition}

\begin{proof}
Since $\mathcal B(\fun)$ is finite, the direct sum of \cref{Eq:Zigzag_Decomp} is in fact a direct product.  Thus, since small limits and categorical products commute \cite[Theorem 3.8.1]{riehl2017category}, \cref{eq:limits_products_commute} holds.  The formula for $\lim \kbb^{I}$ follows from \cref{Eq:Concrete_Limit_For_Finite_Poset}.
\end{proof}

\begin{example}
Let $P$ be the zigzag poset from \cref{ex:Zigzag}, whose Hasse diagram is
\begin{center}
\begin{tikzcd}[ampersand replacement=\&,row sep=1.7ex,column sep=2.5ex,scale=.3,severy label/.append style={font=\small}]
\& v \arrow{r}  \& w                      \&                   \&                                                                                               \\
\&   \& x \arrow[]{r}\arrow[]{u}   \& y 	 					            \\
\&    \&                                                             \& z.\arrow[]{u}
\end{tikzcd}
\end{center}
Note that $\{v,w,x\}$ is a downset of $P$, but $\{v,w,x,y\}$ is not.  Therefore, the following two diagrams have limits $\kbb$ and 0, respectively:
\begin{center}
\begin{tikzcd}[ampersand replacement=\&,row sep=1.7ex,column sep=2.5ex,scale=.3,severy label/.append style={font=\small}]
\& \kbb \arrow{r}  \& \kbb                      \&                   \&                                                                                               \\
\&   \& \kbb \arrow[]{r}\arrow[]{u}   \& 0 	 					            \\
\&    \&                                                             \& 0\arrow[]{u}
\end{tikzcd}
\hskip10pt
\begin{tikzcd}[ampersand replacement=\&,row sep=1.7ex,column sep=2.5ex,scale=.3,severy label/.append style={font=\small}]
\& \kbb \arrow{r}  \& \kbb                      \&                   \&                                                                                               \\
\&   \& \kbb \arrow[]{r}\arrow[]{u}   \& \kbb 	 					            \\
\&    \&                                                             \& 0\arrow[]{u}
\end{tikzcd}
\end{center}
\end{example}

By \cref{Prop:Limit_of_Zigzag}, to compute $\lim \fun$ in the setting of \cref{Thm:Zigzag_decomposition}, it suffices to compute an interval decomposition of $\fun$.  However, note that to compute a limit cone for $\fun$ (and not only for some functor naturally isomorphic to $\fun$), it is not enough to compute $\mathcal{B}(\fun)$; we need an explicit choice of an interval decomposition.

\begin{proof}[Proof of \cref{Thm:Chain_complex_limit}\,(i)]
By \cref{Rem:Zigzag_N^2}, $Q$ has a unique initial scaffold $j\colon P\hookrightarrow Q$, where $P$ is a finite zigzag poset.
Let $\fun$ be the restriction of $H$ to $P$, i.e., $\fun=H \circ i\circ j$.  By \cref{Thm:Computing_Hull_Skeleton_Intervals}\,(i), we can compute $P$ in time $O(n\log n)$.

Let $([f],[g])$ denote the $(\N^2,r)$-complex given as input, and recall that this represents a chain complex of free  $\N^2$-modules.  We compute the restriction of this chain complex to $P$, representing this in memory as a pair of labeled matrices $(\Phi,\Gamma)$, as follows: For a column or row of either $[f]$ or $[g]$ with label $l\in \N^2$, let $l[P]=\Up{l,\N^2}\cap P$; either $l[P]$ is empty or it is both an interval in $P$ and an upset of $P$.  We take the unlabeled matrix underlying $\Phi$ to be the submatrix of $[f]$ formed by the rows and columns of $[f]$ whose labels $l$ satisfy $l[P]\ne\emptyset$.  We take the label of each such row or column of $\Phi$ to be the pair of endpoints of $l[P]$, where the two endpoints are both $p$ when $l[P]=\{p\}$ is a singleton.  The matrix  $\Gamma$ is given as a submatrix of $[g]$ in the analogous way, also defining the labels of $\Gamma$ analogously.  Given $P$, computing the pair $(\Phi,\Gamma)$ takes time $O(r\log n)$ using binary search, assuming we do not copy any rows or columns of $[f]$ or $[g]$.

Given $P$ and $(\Phi,\Gamma)$, the zigzag persistence algorithm of Dey, Hou, and Morozov \cite{DTM25} computes a set of chain representatives of an interval decomposition of $\fun$, with the chain representative of each interval $I\in \mathcal B(\fun)$ specified explicitly at each $p\in I$.  This computation  takes time $O(r^3+|P|r^2)= O(r^3+nr^2)$.   (As stated, the algorithm of  \cite{DTM25}  assumes the input to be a zigzag filtration of a simplicial complex, but it extends immediately to a chain complex represented in this format.) \cref{Prop:Limit_of_Zigzag} then yields a presection basis for $\lim \fun$ in time $O(nr^2)$.  Hence, this approach computes $\lim \fun$, and thus $\lim(H\circ i)$, in total time $O(n\log n+r^3+nr^2)$. \end{proof}

\begin{remark}\label{Rem:Drop_nr2_term}
When stored in the naive way, the presection basis output by the computation of \cref{Thm:Chain_complex_limit}\,(i) requires $\Theta(nr^2)$ space in the worst case.  However, if we output the presection basis in a simple compressed format, then only $O(n+r^3)$ space is needed, and with a minor modification to the algorithm given above, the term $nr^2$ in the bound of \cref{Thm:Chain_complex_limit}\,(i) can be removed.

In brief, this works as follows: Before running the zigzag persistence algorithm of \cite{DTM25}, we quotient the poset $P$ to a zigzag poset $P'$ of size $O(r)$ by collapsing to a single point
each subposet $\bullet \rightarrow \bullet \leftarrow \bullet$ for which the structure maps in the chain complex are identities.  Let $\rho\colon \Ob(P)\to \Ob(P')$ be the quotient map on objects.  By composing structure maps and their inverses, the chain complex of $P$-modules represented by $(\Phi,\Gamma)$ descends to a chain complex of $P'$-modules; this chain complex is represented by the pair of labeled matrices obtained by replacing each label $(p,q)$ of $(\Phi,\Gamma)$ with $(\rho(p),\rho(q))$.  We denote the homology of the latter chain complex as $\fun'$.
Given $P'$ and $(\Phi,\Gamma)$, the algorithm of \cite{DTM25} computes a set $B'$ of chain representatives of an interval decomposition of $\fun'$ in time $O(r^3+|P'|r^2)=O(r^3)$.  The set $B'$ pulls back under $\rho$ to a set $B$ of chain representatives of an interval decomposition of $\fun$.
To represent $B$ in memory, it suffices to store $B'$ and $\rho$.
The presection basis obtained from $B$ via \cref{Prop:Limit_of_Zigzag} then admits an analogous representation, computed by taking a subset of $B'$.  If we output the presection basis in this compressed format, then the total cost of computing $\lim \fun$ is $O(n\log n+r^3)$.
\end{remark}

\begin{remark}\label{Rem:No_Limit_in_MM_Time}
If we wish to only compute $\dim (\lim (H\circ i))$ or the limit cone of a functor naturally isomorphic to $H\circ i$, then the complexity bound of \cref{Thm:Chain_complex_limit}\,(i)  can be improved to $O(n\log n+r^{\omega})$ by using either the zigzag persistence algorithm of Milosavljević, Morozov, and Skraba \cite{milosavljevic2011zigzag} or the more practical algorithm of Dey and Hou \cite{DH22}.  However, these algorithms do not provide an explicit interval decomposition, which we need to compute the cone maps of $\lim (H\circ i)$.
\end{remark}

\section{Proofs of Results on Generalized Rank Computation}\label{Sec:Main_Compute_Gen_Rank}
In this section, we prove \cref{Thm:Grank_Gen_Hull_Skeleta} and \cref{Thm:Grank_Nd_Hull_Skeleta}, our bounds on the cost of generalized rank computation.  We use the following notation:

\begin{notation}\label{Notation:Sec_8_h_notation}
Given a functor $\fun\colon Q\to \Vec$ where $Q$ is a finite, connected poset, let $j^{\mathcal I}\colon P^{\mathcal I}\hookrightarrow Q$ and $j^{\mathcal F}\colon P^{\mathcal F}\hookrightarrow Q$ be initial and final scaffolds of $Q$.  By \cref{Cor:Initial_Scaffold_Connected} and its dual, both $P^{\mathcal I}$ and $P^{\mathcal F}$ are connected.  Given a choice of a minimal element $m\in P^{\mathcal I}$ and a maximal element $w\in P^{\mathcal F}$ with $m\leq w$ in $Q$, let $\{m\leq w\}$ denote the two-element poset with a single non-identity relation $m\leq w$.  Consider the connected poset
 \[P\coloneq P^{\mathcal I}\sqcup P^{\mathcal F} \cup \{m\leq w\},\]
and let $j\colon P\rightarrow Q$ be the obvious functor; that is, $j$ restricts to $j^{\mathcal I}$ and $j^{\mathcal F}$ on $P^{\mathcal I}$ and $P^{\mathcal F}$, respectively, and $j(m\leq w)=(m\leq w)$.
\end{notation}

Note that the functor $j$ defined above needn't be injective on objects, as $P^{\mathcal I}\cap P^{\mathcal F}$ can be non-empty.

We first prove \cref{Thm:Grank_Gen_Hull_Skeleta}.  Recall the statement:

\begin{repcorollary}{Thm:Grank_Gen_Hull_Skeleta}
Let $Q$ be a finite, connected poset with $n$ extrema and let $P^{\mathcal I}$, $P^{\mathcal F}$ be initial and final scaffolds of $Q$.  Given the Hasse diagram $(\mathcal V,\mathcal E)$ of $Q$ and a $(Q,r)$-complex with homology $H$, we can compute $\grank(H)$ in time
\[O(|\mathcal V||\mathcal E|+(|\hom P^{\mathcal I}|+|\hom P^{\mathcal F}|)n^{\omega-1}r^\omega).\]
\end{repcorollary}

\begin{proof}
To compute $\grank(H)$, we first compute initial and final scaffolds $j^{\mathcal I}\colon P^{\mathcal I}\hookrightarrow Q$ and $j^{\mathcal F}\colon P^{\mathcal F}\hookrightarrow Q$ via \cref{Thm:Compute_Hull_Skeleton}.  Choosing $m\leq w$ as in \cref{Notation:Sec_8_h_notation} and letting $j\colon P\rightarrow Q$ be as defined there, we compute a matrix representation of $H\circ j$ via \cref{Prop:Chain_Complex_Compute_Matrix_Rep,Rem:Compute_all_relations}, which, in particular, yields a basis $B_p$ of $H_p$ for each $p\in P$.  This matrix representation restricts to matrix representations of $H\circ j^{\mathcal I}$ and $H\circ j^{\mathcal F}$, which we use to compute $\lim(H)$ and $\colim(H)$ as we did in \cref{Thm:Chain_complex_limit_gen}.

Let $\delta_m\colon \lim H \to H_m$ and $\delta'_w\colon H_w\to \colim H$ be the respective cone and cocone maps.  Observe that \[\grank(H)=\rank( \delta'_w\circ H_{mw} \circ \delta_m).\]   Since $m$ is a minimal element of $P^{\mathcal I}$, hence a minimal element of $Q$, a matrix representation $[\delta_m]$ of $\delta_m$ with respect to the basis $B_m$ is given by coordinate projection of the computed presection basis of $\lim H$; see \cref{Rems:Def_of_Lim_Computation}\,(i).   Dually, coordinate projection also yields a matrix representation $[\delta'_w]$ of $\delta'_w$ with respect to the basis $B_w$.  The matrix representation of $H\circ j$ furnishes a matrix representation $[H_{mw}]$ of $H_{mw}$ with respect to the bases $B_m$ and $B_w$.  To compute $\grank(H)$, we compute the matrix product $[\delta'_w][H_{mw}][\delta_m]$ and perform Gaussian elimination on this.

The cost of these computations is dominated by the cost of computing $\lim H$ and $\colim H$, which is given by \cref{Thm:Chain_complex_limit_gen}.  The result follows.
\end{proof}

We next turn to the proof of \cref{Thm:Grank_Nd_Hull_Skeleta}.  Recall the statement:

\begin{repcorollary}{Thm:Grank_Nd_Hull_Skeleta}
Given an $(\N^d,r)$-complex with homology $H$ and the $n$ extrema of a finite interval $i\colon Q\hookrightarrow \N^d$, we can compute $\grank(H\circ i)$ in time
\begin{itemize}
\item[(i)] $O(n\log n+r^{\omega})$ for $d=2$,
\item[(ii)] $O((nr)^{\omega})$ for $d=3$,
\item[(iii)] $O(n^4+n^{\omega+1}r^{\omega})$ for $d>3$.
\end{itemize}
\end{repcorollary}

\begin{proof}[Proof of \cref{Thm:Grank_Nd_Hull_Skeleta}\,(ii,iii)]
The proof is the same as the proof of \cref{Thm:Grank_Gen_Hull_Skeleta}, using \cref{Thm:Computing_Hull_Skeleton_Intervals} in place of \cref{Thm:Compute_Hull_Skeleton}, as well as
\cref{Thm:Chain_complex_limit,rem:maxima_instead_of_upset_pres} in place of \cref{Thm:Chain_complex_limit_gen}.
\end{proof}

In the setting of \cref{Thm:Grank_Nd_Hull_Skeleta}\,(i), the proof of \cref{Thm:Grank_Gen_Hull_Skeleta}, together with \cref{Thm:Chain_complex_limit}\,(i) and \cref{Rem:Drop_nr2_term}, yields the bound $O(n\log n+r^3)$ for generalized rank computation over a finite interval in $\N^2$, but does not yield the stronger bound  $O(n\log n+r^{\omega})$ given by \cref{Thm:Grank_Nd_Hull_Skeleta}\,(i).  The term $r^3$ arises from the cost of (co)limit computation via \cref{Thm:Chain_complex_limit}\,(i); see \cref{Rem:No_Limit_in_MM_Time}.
We therefore need a different argument to prove  \cref{Thm:Grank_Nd_Hull_Skeleta}\,(i).  Our proof hinges on the following proposition.

\begin{proposition}\label{Lem:grank_initial_final}
In the setting of \cref{Notation:Sec_8_h_notation}, we have \[\grank(\fun)=\grank (\fun\circ j).\]
\end{proposition}

We will only need the case of~\cref{Lem:grank_initial_final} where $Q$ is an interval in $\N^2$, but it is no more difficult to prove the general result.  As noted in \cref{Sec:Intro}, a variant of the proposition was obtained independently in \cite{brustle2025generalized}.

\begin{proof}[Proof of \cref{Lem:grank_initial_final}]
Let $k^{\mathcal I}:P^{\mathcal I}\hookrightarrow P$ and $k^{\mathcal F}:P^{\mathcal F}\hookrightarrow P$ be the inclusions. Note that $j^{\mathcal I}=j\circ k^{\mathcal I}$ and $j^{\mathcal F}=j\circ k^{\mathcal F}$.  By \cref{theorem:MIS_main}, $j^{\mathcal I}$ is initial, so it induces an isomorphism \[\lim \fun\to \lim(\fun\circ j^{\mathcal I})=
 \lim(\fun\circ j\circ k^{\mathcal I}).\]  By the functoriality of limits with respect to the index category (see \cref{Sec:Limits}), this isomorphism factors as
\[
\begin{tikzcd}
\lim \fun \ar[r, hook] & \lim (\fun\circ j) \ar[r, two heads,"\rho"] & \lim (\fun\circ j^{\mathcal I})
\end{tikzcd}
\]
where the first map is an injection and the second map is a surjection.
Dually, the map $j^{\mathcal F}$ induces an isomorphism
$\colim(\fun\circ j^{\mathcal F})\to \colim \fun$ factoring as
\[
\begin{tikzcd}
\colim (\fun\circ j^{\mathcal F}) \ar[r, hook,"\iota"] & \colim (\fun\circ j) \ar[r, two heads] & \colim \fun
\end{tikzcd}
 \]
 where again, the first map is an injection and the second map is a surjection.
 Let $\alpha$ denote the  composite map
\[\lim (\fun\circ j) \to \fun_m\to \fun_{w} \to \colim (\fun\circ j)\]
and let $\beta$ denote the composite map
\[\lim( \fun\circ j^{\mathcal I}) \to \fun_m \to \fun_{w}\to \colim (\fun \circ j^{\mathcal F}). \]
Then the following diagram commutes:
\[ \begin{tikzcd}
\lim(\fun \circ j) \arrow[dr] \arrow[dd, two heads,"\rho"] \arrow[rrr, "\alpha"'] &
&
&
\colim(\fun \circ j) \\
&
\fun_m \arrow[r] &
\fun_w \arrow[ur] \arrow[dr] &
\\
\lim(\fun \circ j^{\mathcal I}) \arrow[ur] \arrow[rrr, "\beta"] &
&
&
\colim(\fun \circ j^{\mathcal F}) \arrow[uu, hook,"\iota"]
\end{tikzcd} \]
Thus, $\alpha=\iota\circ \beta\circ\rho$.  Since $\rho$ is a surjection and $\iota$ is an injection, we have $\rank(\alpha)=\rank(\beta)$.  By definition, we have $\grank(\fun \circ j)=\rank(\alpha)$.  Moreover, since $j^{\mathcal I}$ and $j^{\mathcal F}$ are initial and final, respectively, \cref{Prop:Initial_Functors_Preserve_Limits} and its dual imply that $\grank(\fun)=\rank(\beta)$.  Therefore,
\[\grank(\fun \circ j)=\rank(\alpha)=\rank(\beta)=\grank(\fun).\qedhere\]
\end{proof}

\begin{proof}[Proof of \cref{Thm:Grank_Nd_Hull_Skeleta}\,(i)]
We use an argument similar to the proof of \cref{Thm:GRank_N2} in \cite{DKM22}, but framed in our general formalism of initial and final scaffolds.  We compute initial and final scaffolds $P^{\mathcal I}$ and $P^{\mathcal F}$ of $Q$ in time $O(n\log n)$ using \cref{Thm:Computing_Hull_Skeleton_Intervals}\,(i).
Recall from \cref{Rem:Zigzag_N^2} that $P^{\mathcal I}$ is a zigzag poset; by duality, the same is true for $P^{\mathcal F}$.  Let $m$ be the leftmost minimal element of $Q$ and let $w$ be the leftmost maximal element of $Q$.  Then $m\in P^{\mathcal I}$, $w\in P^{\mathcal F}$, and $m\leq w$.  Let $P$ and $j$ be as in \cref{Notation:Sec_8_h_notation}, and let $\fun=H\circ i$.  Then $P$ is a zigzag poset and the barcode of the zigzag module $\fun \circ j$ can be computed in time $O(n\log n+r^\omega)$ via the algorithm of either \cite{milosavljevic2011zigzag} or \cite{DH22}.  The number of copies of the interval $P$ in the barcode is the number of interval summands of $\fun\circ j$ with support $P$, which in turn is $\grank(\fun\circ j)$ by \cite[Corollary 3.2]{chambers2018persistent}.
We have $\grank(\fun\circ j)=\grank(\fun)=\grank(H\circ i)$ by \cref{Lem:grank_initial_final}.  The result follows.
\end{proof}

\section{Questions for Future Work}\label{Sec:Future_Work}
Our main results concern initial functors with target a poset.  It would be interesting to extend these results to initial functors with target an arbitrary finite category.  It would also be worthwhile to study the optimality of our algorithmic results, particularly for limit computation.  The algorithms of this paper for computing limits and generalized ranks should be relatively easy to implement using ordinary (cubic-time) Gaussian elimination; a good implementation would open the door to applications.  In particular, it would be interesting to explore applications to the featurization of 3-parameter persistent homology for supervised learning, along the lines of what Xin et al. \cite{xin2023gril} have done in the 2-parameter case.

\bibliographystyle{abbrv}
\bibliography{biblio}
\end{document}